\documentclass[10pt]{amsart}

\usepackage{eucal,url,amssymb,stmaryrd,enumerate,amscd,}
\usepackage[pagebackref,colorlinks=true]{hyperref}
\usepackage{amsfonts}
\usepackage{amsmath,amsthm,amssymb,amscd,enumerate,eucal,url,stmaryrd,verbatim}

\numberwithin{equation}{section}
\newtheorem{thrm}{Theorem}[section]
\newtheorem{lemma}[thrm]{Lemma}
\newtheorem{prop}[thrm]{Proposition}
\newtheorem{cor}[thrm]{Corollary}
\newtheorem{dfn}[thrm]{Definition}

\newtheorem{rmrk}[thrm]{Remark}

\newtheorem{conv}[thrm]{Convention}
\newcommand{\Hn}{\mathbb{H}^n}

\newcommand{\QH}{\boldsymbol {G\,(\mathbb{H})}}

\newcommand{\abs}[1]{\lvert #1 \rvert}

\title[Conformal quaternionic contact curvature and the local sphere theorem]
{Conformal quaternionic contact curvature and the local sphere theorem.\\
La courbure conforme d'une structure de contact quaternionienne et structures localment plates.}

\begin{document}

\begin{abstract}
A tensor invariant is defined on a quaternionic contact  manifold in terms of the
curvature and torsion of the Biquard connection involving  derivatives up to
 third order of the contact form. This tensor,
called quaternionic contact conformal curvature, is similar to
the Weyl conformal curvature in Riemannian geometry and to the
Chern-Moser tensor in CR geometry.
It is shown that a quaternionic contact
manifold is locally quaternionic contact conformal to the standard flat quaternionic contact structure
on the quaternionic Heisenberg group, or equivalently, to the
standard 3-sasakian structure on the sphere iff the quaternionic contact
conformal curvature vanishes.

 Un tenseur est d\'{e}fini sur une vari\'{e}t\'{e} avec une structure de contact quaternionienne  en utilisant la courbure et la torsion de la connexion de Biquard. Ce tenseur, appel\'{e}e la courbure conforme d'une structure de contact quaternionienne, ne dependent que des d\'{e}riv\'{e}s de la  troisi\`{e}me ordre de  form de contact et qui est similaire \`{a}  la courbure de Weyl dans  le cas riemannienne et \`{a} le tenseur de Chern-Moser  dans la g\'{e}om\'{e}trie CR. Il est d\'{e}montr\'{e} que une structure de contact quaternionienne est localement conforme \`{a} la structure de contact quaternionienne plate sur le groupe de Heisenberg, ou encore, \`{a} la  structure 3-sasakienne sur la sph\`{e}re quaternionic si et seulement si  la courbure conforme de contact quaternionienne est nulle.

\end{abstract}


\date{\today}
\thanks{This project has been funded in part by the National Academy of
Sciences under the [Collaboration in Basic Science and Engineering Program 1
Twinning Program] supported by Contract No. INT-0002341 from the National
Science Foundation. The contents of this publication do not necessarily
reflect the views or policies of the National Academy of Sciences or the
National Science Foundation, nor does mention of trade names, commercial
products or organizations imply endorsement by the National Academy of
Sciences or the National Science Foundation.}
\author{Stefan Ivanov}
\address[Stefan Ivanov]{University of Sofia, Faculty of Mathematics and
Informatics, blvd. James Bourchier 5, 1164, Sofia, Bulgaria}
%\address{and Max-Planck-Institut f\"ur Mathematik\\
%Vivatsgasse 7\\
%D-53111 Bonn\\Germany}
\email{ivanovsp@fmi.uni-sofia.bg}

\author{Dimiter Vassilev}
\address[Dimiter Vassilev]{ Department of Mathematics and Statistics\\ University of New Mexico\\
Albuquerque, New Mexico, 87131}
%\address{University of California, Riverside\\
%Riverside, CA 92521}
%\address{Max-Planck-Institut f\"ur Mathematik\\
%Vivatsgasse 7\\
%D-53111 Bonn\\
%Germany}
\email{vassilev@math.unm.edu}

\maketitle

\tableofcontents

\setcounter{tocdepth}{2}

\section{Introduction}

{It is well known that the sphere at infinity of a  non-compact
symmetric space $M$ of rank one carries a natural
Carnot-Carath\'eodory structure, see \cite{M,P}.} { A quaternionic
contact (qc) structure, introduced in \cite{Biq1,Biq2}, appears
naturally as the conformal boundary at infinity of the
quaternionic hyperbolic space. Such structures have been
considered in connection with the quaternionic contact Yamabe
problem, \cite{GV,Wei,IMV,IMV1}. A particular case of this problem
amounts to find the extremals and the best constant in the $L^2$
Folland-Stein Sobolev-type embedding, \cite{F2} and \cite{FS}. A
complete description of the extremals and the best constant  on
the seven dimensional quaternionic Heisenberg group was given in
\cite{IMV1}.}

{A qc structure on a real (4n+3)-dimensional manifold $M$ is a codimension
three distribution $H$ locally given as the kernel of 1-form $%
\eta=(\eta_1,\eta_2,\eta_3)$ with values in $\mathbb{R}^3$ and the three
2-forms $d\eta_i|_H$ are the fundamental 2-forms of a quaternionic structure
on $H$, i.e., there exists a Riemannian metric $g$ on $H$ and three local
almost complex structures $I_i$ on $H$ satisfying the commutation relations
of the imaginary quaternions, $I_1I_2I_3=-1$, such that, $%
d\eta_i|_H=2g(I_i.,.)$ . The 1-form $\eta$ is determined up to a conformal
factor and the action of $SO(3)$ on $\mathbb{R}^3$, and therefore $H$ is
equipped with a conformal class $[g]$ of Riemannian metrics and a 2-sphere
bundle of almost complex structures, the quaternionic bundle $\mathbb{Q}$.
The 2-sphere bundle of one forms determines uniquely the associated metric
and a conformal change of the metric is equivalent to a conformal change of
the one forms. To every metric in the fixed conformal class one can
associate a linear connection $\nabla$ preserving the qc structure, see \cite{Biq1},
which we shall call the Biquard connection. }

The transformations preserving a given quaternionic contact
structure $\eta$, i.e. $\bar\eta=\mu\Psi\eta$ for a positive smooth
function $\mu$ and a $SO(3)$ matrix $\Psi$ with smooth functions as
entries, are called \emph{quaternionic contact conformal (qc
conformal for short) transformations}.

Examples of QC manifolds are given in \cite{Biq1,Biq2,IMV,D1,AFIV}. In
particular, any totally umbilic hypersurface of a quaternionic
K\"ahler or hyper K\"ahler manifold carries such a structure
\cite{IMV}. A basic example
is provided by any 3-Sasakian manifold which can be defined as a $(4n+3)$%
-dimensional Riemannian manifold whose Riemannian cone is a hyper
K\"ahler manifold.  It was shown in \cite{IMV} that when the horizontal
scalar curvature $Scal$ of the Biquard connection { (qc scalar curvature for short)} is positive the
torsion endomorphism of the Biquard connection is the obstruction
for a given qc-structure to be locally 3-Sasakian.

The quaternionic Heisenberg group $\QH$ with its "standard"
left-invariant qc structure is the unique (up to a $SO(3)$-action)
example of a qc structure with flat Biquard connection \cite{IMV}.
{As a manifold $\QH \ =\mathbb H^n\times\text {Im}\, \mathbb{H}$,
while the group multiplication is given by
\[
( q', \omega')\ =\ (q_o, \omega_o)\circ(q, \omega)\ =\ (q_o\ +\ q, \omega\ +\ \omega_o\ +
\ 2\ \text {Im}\  q_o\, \bar q),
\]
\noindent where $q,\ q_o\in\mathbb H^n$ and $\omega, \omega_o\in \text
{Im}\, \mathbb{H}$. The standard flat quaternionic contact structure
is defined by the left-invariant quaternionic contact form
$\tilde\Theta\ =\ (\tilde\Theta_1,\ \tilde\Theta_2, \
\tilde\Theta_3)\ =\ \frac 12\ (d\omega \ - \ q' \cdot d\bar q' \ + \
dq'\, \cdot\bar q')$, where $.$ denotes the quaternion
multiplication.}

The aim of this paper is to find a tensor invariant on the tangent
bundle characterizing locally the qc structures which are
quaternionic contact conformally equivalent to the flat
qc-structure. With this goal in mind, we describe a curvature-type
tensor $W^{qc}$ defined in terms of the curvature and torsion of
the Biquard connection by \eqref{qccm} { involving  derivatives up to
second order of the horizontal metric}, whose form is
similar to the Weyl conformal curvature in Riemannian geometry
(see e.g. \cite{Eis}) and to the Chern-Moser invariant in CR
geometry \cite{ChM}, see also \cite{W}. We call $W^{qc}$ the
\emph{quaternionic contact conformal curvature, qc conformal
curvature for short}. The main purpose of this article is to prove
the following two facts.
\begin{thrm}\label{main1}
The qc conformal curvature $W^{qc}$ is invariant under quaternionic contact conformal transformations.
\end{thrm}
\begin{thrm}\label{main2}
A qc structure on a (4n+3)-dimensional smooth manifold is locally quaternionic contact conformal to the
standard flat qc structure on the quaternionic Heisenberg group $\QH$
if and only if the qc conformal curvature vanishes,
$W^{qc}=0$.
\end{thrm}

The quaternionic Cayley transform establishes a conformal
quaternionic contact { isomorphism} between the standard 3-Sasaki
structure on the { punctured} sphere $S^{4n+3}$ and the standard
qc structure on $\QH$ \cite{IMV}, { which combined with %As ansequence
Theorem~\ref{main2}  gives }%we obtain}
\begin{cor}\label{main3}
A QC manifold is locally quaternionic contact conformal to the standard 3-Sasakian 
sphere $S^{4n+3}$ if and only if the qc conformal curvature
vanishes, $W^{qc}=0$.
\end{cor}
We note that for locally 3-Sasakian manifolds a curvature invariant under  very special quaternionic contact
conformal transformations, which  preserve the 3-Sasakian condition, is defined in \cite{AK}. It is shown that
the vanishing of this invariant is equivalent to the structure being locally isometric to the 3-Sasaki
structure on the sphere. In particular, this shows that the standard 3-Sasakian structure on the sphere is
locally rigid with respect to qc conformal transformations preserving the 3-Sasakian condition.

We consider the question of local flatness in its full generality following  the classical approach used by
H.Weyl, see e.g. \cite{Eis}, while \cite{ChM}, \cite{W} and \cite{AK} followed the Cartan method of
equivalence. Indeed, Theorem~\ref{main1} and Theorem~\ref{main2} could be deduced following the Cartan method
of equivalence as in \cite{ChM}. Vice versa, the Chern-Moser tensor \cite{ChM} can be determined (c.f.
\cite{IVZ}) as an obstruction to the pseudohermitian flatness following the approach used in the present paper.

\begin{rmrk}
Following the work of Cartan and Tanaka, a qc structure can be considered as an example of what has become known as a parabolic geometry.
The quaternionic Heisenberg group, as well as, { the (4n+3) dimensional sphere}, due to the Cayley transform, provide the flat models of
such a geometry. It is well known that the curvature of the corresponding {regular} Cartan connection is the obstruction for the local
flatness. However, the Cartan curvature is not a tensor field on the tangent bundle and it is highly nontrivial to extract a tensor field
involving the lowest order derivatives of the structure which implies the vanishing of the obstruction. { Theorem 1.2 suggests that a
necessary and sufficient condition for the vanishing of the Cartan curvature of a qc structure is the vanishing of the qc-conformal
curvature tensor, $W^{qc}=0$.}
\end{rmrk}

In the concluding section of the paper we apply our results in a
standard way to show a Ferrand-Obata type theorem concerning the
conformal quaternionic contact automorphism group. Such result was
proved in the general context of parabolic geometries admitting {
regular} Cartan connection in \cite{F}.

It is  expected that the conformal quaternionic contact curvature
tensor will be a useful tool in the analysis of the quaternionic
contact Yamabe problem, see \cite{Biq1},  \cite{Wei}, \cite{IMV}
and \cite{IMV1}. According to \cite{Wei} the qc Yamabe constant of
a compact qc manifold is less or equal than that of the sphere.
Furthermore, if the constant is strictly less than that of the
sphere the qc Yamabe problem has a solution, i.e., there is a
{ global qc conformal transformation sending the given qc structure to a qc structure  with constant
qc scalar curvature}. A natural conjecture is that the qc Yamabe constant of
every compact locally non-flat manifold (in conformal quaternionic
contact sense) is strictly less than the qc Yamabe constant of the
sphere with its standard qc structure. Recall that the qc Yamabe
constant of $(M, [\eta])$ is
$$ \lambda(M, [\eta])\ =\ \inf \{ \Upsilon (u)\ =\
\int_M\Bigl(4\frac {Q+2}{Q-2}\ \abs{\nabla u}^2\ +\ \text{Scal}\,
u^2\Bigr) dv_g :\ \int_M u^{2^*}\, dv_g \ =\ 1, \ u>0\}.$$ Here
$dv_g$ denotes the Riemannian volume form of the Riemannian metric
on $M$ obtained by extending in a natural way the horizontal
metric associated to $\eta$ and $|.|$ is the horizontal norm.
Guided by the conformal and CR cases the tensor $W^{qc}$ should be
instrumental for the proof of the above conjecture.
\begin{rmrk}
{After the paper was posted on the arXiv, Kunkel \cite{Ku}
constructed quaternionic contact parabolic normal coordinates and
showed that the only non-trivial scalar invariant of weight at
most four is the square of the horizontal norm of the qc conformal
curvature. This fact suggests that a suitable asymptotic expansion
of the qc Yamabe functional could give an expression of the qc
Yamabe constant in terms of the qc Yamabe constant of the sphere
and the horizontal norm of the qc conformal curvature  leading to a
proof of the above conjecture.}
\end{rmrk}

{ Two examples where Theorem~\ref{main2} applies are given in
\cite{AFIV}.
\begin{comment}Recall that the two-step nilpotent seven
dimensional quaternionic Heisenberg Lie algebra is defined by the
following structure equations for the left invariant 1-forms
\begin{equation*}%\label{7heis}
\begin{aligned}
de^1 = de^2 =  de^3 = de^4 = 0,\quad de^5 = 2 e^{12} + 2 e^{34},
\quad de^6 = 2 e^{13} - 2 e^{24},\quad de^7 = 2 e^{14} +2 e^{23},
\end{aligned}
\end{equation*}
where, as usual, $e^{ij}$ denotes $de^i\wedge de^j$. \end{comment}
%In the first example we

Consider the simply connected Lie group
$G_1$ with structure equations %for the left invariant 1-forms given
%by
%\begin{equation*}%\label{ex11}
\begin{align*}
& de^1=0\qquad de^2=-e^{12}-2e^{34}-\frac12e^{37}+\frac12e^{46}\\ &
de^3=-e^{13}+2e^{24}+\frac12e^{27}-\frac12e^{45}\qquad
de^4=-e^{14}-2e^{23}-\frac12e^{26}+\frac12e^{35} \\ &
de^5=2e^{12}+2e^{34}-\frac12e^{67},
\qquad de^6=2e^{13}+2e^{42}+\frac12e^{57}\qquad
 de^7=2e^{14}+2e^{23}-\frac12e^{56}
\end{align*}
%\end{equation*}
where, as usual, $e^{ij}$ denotes $e^i\wedge e^j$. It is shown in
\cite{AFIV} that $H=span\{e^1,\dots, e^4\}$,
$\eta_1=e^5, \quad \eta_2=e^6, \quad \eta_3=e^7, \quad
 \omega_1=e^{12}+e^{34}, \quad
\omega_2=e^{13}+e^{42}, \quad \omega_3=e^{14}+e^{23}$
\noindent determine  a (global) qc structure on $G_1$, for which the
torsion endomorphism of the Biquard connection vanishes, but the
qc-structure is not locally 3-Sasakian since the qc scalar curvature
is negative. Furthermore, the qc conformal curvature vanishes, $W^{qc}=0$, hence
due to Theorem~\ref{main2}, the considered qc structure on $G_1$  is
locally qc conformally equivalent to the flat qc structure on the
7-dimensional quaternionic Heisenberg group.

Another example is provided by the simply connected Lie group $G_3$
with structure equations
%\begin{equation*}%\label{ex31}
\begin{align*}
& de^1=-\frac32e^{13}+\frac32e^{24}-\frac34e^{25}+\frac14e^{36}-\frac14e^{47}+\frac18e^{57}\\
&
de^2=-\frac32e^{14}-\frac32e^{23}+\frac34e^{15}+\frac14e^{37}+\frac14e^{46}-\frac18e^{56}\qquad
de^3=0 \\
& de^4=e^{12}+e^{34}+\frac12e^{17}-\frac12e^{26}+\frac14e^{67}
\qquad
de^5=2e^{12}+2e^{34}+e^{17}-e^{26}+\frac12e^{67} \\
& de^6=2e^{13}+2e^{42}+e^{25},\qquad de^7=2e^{14}+2e^{23}-e^{15}.
\end{align*}
%\end{equation*}
A global qc structure is determined by $H=span\{e^1,\dots, e^4\}$,
$\eta_1=e^5, \quad \eta_2=e^6, \quad \eta_3=e^7, \quad
 \omega_1=e^{12}+e^{34}, \quad
\omega_2=e^{13}+e^{42}, \quad \omega_3=e^{14}+e^{23}$.
In this case, by \cite{AFIV}, the
torsion endomorphism and the qc conformal curvature do not vanish.
In particular, this qc structure on $G_3$ is not locally qc
conformally flat according to Theorem~\ref{main2}.}

\textbf{Organization of the paper.} The paper relies heavily on the Biquard
connection introduced in
\cite{Biq1} and the properties of its torsion and curvature discovered in \cite{IMV}.
In order to make the present paper self-contained, in Section \ref{s:review}
we give a review of the notion of a quaternionic contact structure and
collect formulas and results from \cite{Biq1} and \cite{IMV} that will be used in the
subsequent sections.
\begin{conv}\label{conven}
We use the following conventions:
\begin{enumerate}[a)]
\item We shall use $X,Y,Z,U$ to denote horizontal vector fields, i.e. $X,Y,Z,U\in H$;
\item $\{e_1,\dots,e_{4n}\}$ denotes an orthonormal basis of the horizontal
space $H$;
\item The summation convention over repeated vectors from the basis $%
\{e_1,\dots,e_{4n}\}$ will be used. For example, for a (0,4)-tensor $P$, the
formula $k=P(e_b,e_a,e_a,e_b)$ means
%\begin{equation*}\label{e:horizontal riccic}
$%\begin{equation*}
k=\sum_{a,b=1}^{4n}P(e_b,e_a,e_a,e_b);
$%\end{equation*}%
%\end{equation*}
\item The triple $(i,j,k)$ denotes any cyclic permutation of
$(1,2,3)$. In particular, any equation involving $i,j,k$ holds for
any such permutation.
\item $s$ and $t$ will be any numbers from the set $\{1,2,3\}, \quad
s,t\in\{1,2,3\}$.
%\item We shall denote with $\nabla h$ the horizontal gradient of the function $h$,
%    while $dh$ means as usual the differential of the function $h$, or the full gradient of $h$.
\end{enumerate}
\end{conv}

\textbf{Acknowledgements} The research was done during the visit of
S.Ivanov in the Max-Plank-Institut f\"ur Mathematics, Bonn and the
final draft of the paper was prepared with both authors residing at
MPIM, Bonn.  The authors thank  MPIM, Bonn for providing the support
and an excellent research environment. S.I. is a Senior Associate to
the Abdus Salam ICTP. The authors  also like to acknowledge The
National Academies for the financial support, and University of
Sofia and University of California Riverside for hosting the respective visits of
the authors which contributed in the writing of the paper.
S.I. is partially supported by the Contract 082/2009 with the University of Sofia `St.Kl.Ohridski` and
Contract "Idei", DO 02-257/18.12.2008. %\vskip.3truein

\section{Quaternionic contact manifolds}

\label{s:review} In this section we will briefly review the basic notions of
quaternionic contact geometry and recall some results from \cite{Biq1} and
\cite{IMV}.

For the purposes of this paper, a quaternionic contact (QC) manifold $(M, g,
\mathbb{Q})$ is a $4n+3$ dimensional manifold $M$ with a codimension three
distribution $H$ equipped with a metric $g$ and an Sp(n)Sp(1) structure,
i.e., we have
\begin{enumerate}
\item[i)] a 2-sphere bundle $\mathbb{Q}$ over $M$ of almost complex structures
$I_s\,:H \rightarrow H,\quad I_s^2\ =\ -1$, satisfying the commutation relations of the imaginary
quaternions $I_1I_2=-I_2I_1=I_3$ and
 $\mathbb{Q}= \{aI_1+bI_2+cI_3:\ a^2+b^2+c^2=1 \}$;
\item[ii)] $H$ is locally the kernel of a 1-form $\eta=(\eta_1,\eta_2,\eta_3)$ with
values in $\mathbb{R}^3$ satisfying the compatibility condition  \hspace{3mm}
$%\begin{equation*}  \label{con1}
2g(I_sX,Y)\ =\ d\eta_s(X,Y).% \quad s=1,2,3, \quad X,Y\in H.
$%\end{equation*}
\end{enumerate}

A QC manifold $(M, \bar g,\mathbb{Q} )$ is called quaternionic contact conformal (
qc-conformal for short) to $(M, g,\mathbb{Q} )$ if $\bar g\in [g]$. In that case, if $\bar\eta$ is a
corresponding associated one-form with complex structures $\bar I_s$,
we have $\bar\eta\ =\ \mu\, \Psi\,\eta$ for some $\Psi\in SO(3)$
and a positive function $\mu$. In particular, starting with a QC manifold $%
(M, \eta)$ and defining $\bar\eta\ =\ \mu\, \eta$ we obtain a QC manifold $%
(M, \bar\eta)$ qc-conformal to the original one.

On a quaternionic contact manifold there exists a canonical connection
defined in \cite{Biq1} when the dimension $(4n+3)>7$, and in \cite{D} in the
7-dimensional case.

\begin{thrm}
\cite{Biq1}\label{biqcon} {Let $(M, g,\mathbb{Q})$ be a quaternionic contact
manifold} of dimension $4n+3>7$ and a fixed metric $g$ on $H$ in the
conformal class $[g]$. Then there exists a unique connection $\nabla$ with
torsion $T$ on $M^{4n+3}$ and a unique supplementary subspace $V$ to $H$ in $%
TM$, such that:
\begin{enumerate}
\item[i)] $\nabla$ preserves the decomposition $H\oplus V$ and the $Sp(n)Sp(1)$-structure
on $H$;%, $\nabla g\
%= \ 0, \nabla \mathbb{Q}\subset\mathbb{Q}$;
\item[ii)] for $X,Y\in H$, one has $T(X,Y)=-[X,Y]_{|V}$;
%\item[iii)] $\nabla$ preserves the $Sp(n)Sp(1)$-structure on $H$, i.e., $\nabla g\
%=\ 0$ and $\nabla \mathbb{Q}\subset\mathbb{Q}$;
\item[iii)] for $\xi\in V$, the endomorphism $T(\xi,.)_{|H}$ of $H$ lies in $%
(sp(n)\oplus sp(1))^{\bot}\subset gl(4n)$;
%\item[v)] the connection on $V$ is induced by the natural identification $\varphi
%$ of $V$ with the subspace $sp(1)$ of the endomorphisms of $H$, i.e. $%
%\nabla\varphi=0$.
\end{enumerate}
\end{thrm}
We shall call the above connection \emph{the Biquard connection}.
%Its torsion
%endomorphism $T(\xi,.)_{|H}$ evaluated on $H$ will be called \emph{the torsion of the
%given qc-structure}.
Biquard \cite{Biq1} also described the supplementary subspace $V$, namely,
{locally }$V$ is generated by vector fields $%
\{\xi_1,\xi_2,\xi_3\}$, such that
\begin{equation}  \label{bi1}
\begin{aligned} \eta_s(\xi_k)=\delta_{sk}, \qquad (\xi_s\lrcorner
d\eta_s)_{|H}=0,\quad (\xi_s\lrcorner d\eta_k)_{|H}=-(\xi_k\lrcorner
d\eta_s)_{|H}. \end{aligned}
\end{equation}
The vector fields $\xi_1,\xi_2,\xi_3$ are called Reeb vector fields or
fundamental vector fields.

{\ If the dimension of $M$ is seven, the conditions \eqref{bi1} do not
always hold. Duchemin shows in \cite{D} that if we assume, in addition, the
existence of Reeb vector fields as in \eqref{bi1}, then Theorem~\ref{biqcon}
holds. Henceforth, by a qc structure in dimension $7$ we shall %always
mean a qc structure satisfying \eqref{bi1}.}

Notice that equations \eqref{bi1} are invariant under the natural $SO(3)$
action. Using the triple of Reeb vector fields we extend $g$ to a metric on $%
M$ by requiring %\hspace{2mm}
%\begin{equation*}\label{e:extend metric}
$span\{\xi_1,\xi_2,\xi_3\}=V\perp H \text{ and } g(\xi_s,\xi_k)=\delta_{sk}.
$ %\end{equation*}
\hspace{2mm} \noindent The extended metric does not depend on the action of $%
SO(3)$ on $V$, but it changes in an obvious manner if $\eta$ is multiplied
by a conformal factor. Clearly, the Biquard connection preserves the
extended metric on $TM, \nabla g=0$. We shall also extend the quaternionic
structure by setting $I_{s|V}=0$. The fundamental 2-forms $\omega_s$
of the quaternionic structure $Q$ are defined by
\begin{equation}  \label{thirteen}
2\omega_{s|H}\ =\ \, d\eta_{s|H},\qquad \xi\lrcorner\omega_s=0,\quad \xi\in
V.
\end{equation}
Due to \eqref{thirteen}, the torsion restricted to $H$ has the form
\begin{equation}  \label{torha}
T(X,Y)=-[X,Y]_{|V}=2\omega_1(X,Y)\xi_1+2\omega_2(X,Y)\xi_2+2\omega_3(X,Y)\xi_3.
\end{equation}
The properties of the Biquard connection are encoded in the properties of
the torsion endomorphism $T_{\xi}=T(\xi,.) : H\rightarrow H, \quad \xi\in V$.
Recall that any endomorphism $\Psi$ of $H$ can be decomposed with respect to the
quaternionic structure $(\mathbb{Q},g)$ uniquely into
$Sp(n)$-invariant
parts as follows \hspace{2mm} %\begin{equation*}\label{New4}
$\Psi=\Psi^{+++}+\Psi^{+--}+\Psi^{-+-}+\Psi^{--+}, $ %\end{equation*}
\hspace{2mm} where $\Psi^{+++}$ commutes with all three $I_i$,
$\Psi^{+--}$ commutes with $I_1$ and anti-commutes with the others
two and etc. \noindent The two $Sp(n)Sp(1)$-invariant components are
given by
\begin{equation}
{\label{New21}} \Psi_{[3]}=\Psi^{+++}, \qquad
\Psi_{[-1]}=\Psi^{+--}+\Psi^{-+-}+\Psi^{--+}.
\end{equation}
\noindent Denoting the corresponding (0,2) tensor via $g$ by the
same letter one sees that the $Sp(n)Sp(1)$-invariant components are
the projections on the eigenspaces of the Casimir operator
\begin{equation}  \label{e:cross}
\dag \ =\ I_1\otimes I_1\ +\ I_2\otimes I_2\ +\ I_3\otimes I_3
\end{equation}
corresponding, respectively, to the eigenvalues $3$ and $-1$, see \cite{CSal}.
If $n=1$ then the space of symmetric endomorphisms commuting with all
$I_i, i=1,2,3$ is 1-dimensional, i.e. the [3]-component of any
symmetric
endomorphism $\Psi$ on $H$ is proportional to the identity, $\Psi_{[3]}=%
\frac{|\Psi|^2}{4}Id_{|H}$.

Decomposing the endomorphism $T_{\xi}\in(sp(n)+sp(1))^{\perp}$ into
symmetric part $T^0_{\xi}$ and skew-symmetric part $b_{\xi}$, $%
T_{\xi}=T^0_{\xi} + b_{\xi} $, we summarize the description of the torsion
due to O. Biquard in the following Proposition. %\medskip
\begin{prop}
\cite{Biq1}\label{torb} The torsion $T_{\xi}$ is completely trace-free,
\begin{equation*}  \label{torb0}
tr\, T_{\xi}=g\,(\,T_{\xi}e_a,e_a)=0, \quad tr\, T_{\xi}\circ
I=g\,(\,T_{\xi}e_a,Ie_a)=0, \quad I\in Q.
\end{equation*}
%where $e_1\dots e_{4n}$ is an orthonormal basis of $H$.
The symmetric part
of the torsion has the properties:
\begin{gather*}%\label{tors1}
%T^0_{\xi_s}I_s=-I_sT^0_{\xi_s},                 \\
\begin{aligned}%\label{tors2}
T^0_{\xi_s}I_s=-I_sT^0_{\xi_s},\quad I_2T_{\xi_2}^{0_{+--}}=I_1T_{\xi_1}^{0_{-+-}},\quad
I_3T_{\xi_3}^{0_{-+-}}=I_2T_{\xi_2}^{0_{--+}},\quad
I_1T_{\xi_1}^{0{--+}}=I_3T_{\xi_3}^{0_{+--}}. \end{aligned}
\end{gather*}
The skew-symmetric part can be represented as follows
$b_{\xi_s}=I_s\textbf{u},$ where $\textbf{u}$ is a traceless
symmetric (1,1)-tensor on $H$ which commutes with $I_1,I_2,I_3$.

If $n=1$ then  $\textbf{u}$ vanishes identically, $ \textbf{u}=0$
and the torsion is a symmetric tensor, $T_{\xi}=T^0_{\xi}$.
\end{prop}
The $sp(1)$-connection 1-forms are defined by
%The covariant derivative of the quaternionic contact structure with respect
%to the Biquard connection and the covariant derivative of the distribution $V
%$ are given by
%\begin{equation*} % \label{der}
$\nabla I_i=-\alpha_j\otimes I_k+\alpha_k\otimes I_j$, or equivalently determined with
$\nabla\xi_i=-\alpha_j\otimes\xi_k+\alpha_k\otimes\xi_j.
$ {The vanishing of the $sp(1)$-connection 1-forms on $H
$ is equivalent to the vanishing of the torsion endomorphism of the Biquard
connection \cite{IMV}.}

\subsection{The Ricci type tensors}

Let $R=[\nabla,\nabla]-\nabla_{[\ ,\ ]}$ be the curvature tensor
of $\nabla$. We shall denote the curvature tensor of type (0,4) by
the same letter, $ R(A,B,C,D):=g(R(A,B)C,D),\  A,B,C,D \in
\Gamma(TM)$. The Ricci tensor and the scalar curvature $Scal$ of
the Biquard connection, called \emph{qc-Ricci tensor} and
\emph{qc-scalar curvature}, respectively, are defined by
\begin{equation*}  \label{e:horizontal ricci}
Ric(X,Y)=R(e_a,X,Y,e_a), \qquad
Scal=Ric(e_a,e_a)=R(e_b,e_a,e_a,e_b).
\end{equation*}
The curvature of the Biquard connection admits also several horizontal
traces, defined in \cite{IMV} by
\begin{gather*}
4n\rho_s(X,Y)=R(X,Y,e_a,I_se_a), \hspace{1mm} 4n\tau_s(X,Y)=R(e_a,I_se_a,X,Y), \hspace{1mm}
4n\zeta_s(X,Y)=R(e_a,X,Y,I_se_a).
\end{gather*}
The $sp(1)$-part of $R$ is determined by the Ricci 2-forms by
\begin{equation}  \label{sp1curv}
R(A,B,\xi_i,\xi_j)=2\rho_k(A,B), \qquad A,B \in \Gamma(TM).
\end{equation}
According to \cite{Biq1} the Ricci tensor restricted to $H$ is a
symmetric tensor. If the trace-free part of the qc-Ricci tensor is
zero we call the quaternionic structure \emph{a qc-Einstein
manifold} \cite{IMV}. It is shown in \cite{IMV} that the trace part
of these Ricci type contractions is proportional to the qc-scalar
curvature and the trace-free part of $\rho_s$, $\tau_s$, $\zeta_s$
vanish for exactly when the manifold is qc-Einstein (see
also Theorem~\ref{sixtyseven} below).

With the help of the operator $\dag$ we introduced in \cite{IMV} two
$Sp(n)Sp(1)$-invariant trace-free symmetric 2-tensors $T^0, U$ on
$H$  as follows
\begin{gather} \label{tor}
T^0(X,Y)\overset{def}{=}g((T_{\xi_1}^{0}I_1+T_{\xi_2}^{0}I_2+T_{%
\xi_3}^{0}I_3)X,Y),\quad U(X,Y)\overset{def}{=}g( \textbf{u}X,Y).
\end{gather}
The tensor $T^0$ belongs to the [-1]-eigenspace while $U$ is in the
[3]-eigenspace of the operator $\dag$, i.e., they have the
properties:
\begin{gather}\label{propt}
T^0(X,Y)+T^0(I_1X,I_1Y)+T^0(I_2X,I_2Y)+T^0(I_3X,I_3Y)=0, \\\label{propu}
3U(X,Y)-U(I_1X,I_1Y)-U(I_2X,I_2Y)-U(I_3X,I_3Y)=0.
\end{gather}
{ As customary, we let $T(A,B,C)=g(T_AB,C)=g(T(A,B),C), \quad A,B,C\in\Gamma(TM)$.}

Applying Proposition~\ref{torb} and equation \eqref{propt}, we
obtain the following

\begin{prop}
\label{tcomp} The symmetric part of torsion
endomorphism of Biquard connection satisfy the relations
\begin{gather}\label{t001}
4g(T^0_{\xi_s}I_sX,Y)=4T^0(\xi_s,I_sX,Y)=T^0(X,Y)-T^0(I_sX,I_sY).
\end{gather}
\end{prop}
It is shown in \cite{IMV} that all horizontal Ricci type
contractions of the curvature of the Biquard connection can be
expressed in terms of the torsion of the Biquard connection. With
slight modification based on Proposition~\ref{tcomp} we collect some
facts from \cite[Theorem~1.3, Theorem~3.12, Corollary~3.14 and
Proposition~4.3]{IMV} in the next

\begin{thrm}
\cite{IMV}\label{sixtyseven} %Let $(M^{4n+3},g,\mathbb{Q})$ be a
%quaternionic contact
On a $(4n+3)$-dimensional QC manifold, $n>1$  %mFor any $X,Y
%in H, s=1,2,3$ the horizontal Ricci tensors and the qc-scalar curvature
%satisfy
the next formulas hold
\begin{equation*}  \label{sixtyfour}
\begin{aligned} Ric(X,Y) \ & =\ (2n+2)T^0(X,Y)
+(4n+10)U(X,Y)+\frac{Scal}{4n}g(X,Y)\\ \rho_s(X,I_sY) \ & =\
-\frac12\Bigl[T^0(X,Y)+T^0(I_sX,I_sY)\Bigr]-2U(X,Y)-%
\frac{Scal}{8n(n+2)}g(X,Y),\\ \tau_s(X,I_sY) \ & =\
-\frac{n+2}{2n}\Bigl[T^0(X,Y)+T^0(I_sX,I_sY)\Bigr]-%
\frac{Scal}{8n(n+2)}g(X,Y),\\ \zeta_s(X,I_sY) \ & =\
\frac{2n+1}{4n}T^0(X,Y)+\frac{1}{4n}T^0(I_sX,I_sY)+\frac{2n+1}{2n}U(X,Y)+
\frac{Scal}{16n(n+2)}g(X,Y)\\ T(\xi_{i},\xi_{j})& =-\frac
{Scal}{8n(n+2)}\xi_{k}-[\xi_{i},\xi_{j}]_{H}\qquad Scal\  =\
-8n(n+2)g(T(\xi_1,\xi_2),\xi_3)\\ T(\xi_i,\xi_j,X) &
=-\rho_k(I_iX,\xi_i)=-\rho_k(I_jX,\xi_j)\\
\rho_{i}(X,\xi_{i})\ & =\ -\frac
{X(Scal)}{32n(n+2))} \ +\ \frac 12\, \left
(-\rho_{i}(\xi_{j},I_{k}X)+\rho_{j}(\xi_{k},I_{i}X)+\rho_{k}(\xi_{i},I_{j}X)
\right). \end{aligned}
\end{equation*}
For $n=1$ the above formulas hold with $U=0$.

In particular, the qc-Einstein condition is equivalent to the
vanishing of the torsion endomorphism of the Biquard connection. If
$Scal>0$ {\ the latter} holds exactly when the qc-structure is
3-sasakian up to a multiplication by a constant and an
$SO(3)$-matrix with smooth entries.
\end{thrm}
We derive from [\cite{IMV}, Theorem~4.8]
\begin{prop}\cite{IMV} On a qc manifold the following formula holds
\begin{equation}  \label{div}
(n-1)(\nabla_{e_a}T^0)(e_a,X)+2(n+2)\nabla_{e_a}U)(e_a,X)-\frac{(n-1)(2n+1)}{8n(n+2)}d(Scal)(X)=0.
\end{equation}
\end{prop}
%\medskip
\subsection{Quaternionic Heisenberg group and the quaternionic Cayley
transform}\label{s:standard structures} Since our goal is to classify
quaternionic contact manifolds locally conformal to the  quaternionic
Heisenberg group we recall briefly its definition together with the
definition of the quaternionic Cayley transform as described in
\cite[Section 5.2]{IMV}.  As a manifold the quaternionic Heisenberg group
of topological dimension $4n+3$ is $\QH \ =\mathbb{H}^n\times\text {Im}\,
\mathbb{H}$. The group law is given by \hspace{2mm} $ ( q', \omega')\ =\
(q_o, \omega_o)\circ(q, \omega)\ =\ (q_o\ +\ q, \omega\ +\ \omega_o\ + \
2\ \text {Im}\ q_o\, \bar q), $ \noindent where $q,\ q_o\in\mathbb{H}^n$
and $\omega, \omega_o\in \text {Im}\, \mathbb{H}$. We can identify the
group  $\QH$  with the boundary $\Sigma$ of a Siegel domain in $\mathbb
H^n\times\mathbb{H}$, $ \Sigma\ =\ \{ (q',p')\in \Hn\times\mathbb{H}\ :\
\Re {\ p'}\ =\ \abs{q'}^2 \}. $ The Siegel domain
 $\Sigma$ carries a natural group structure and the map $(q,
\omega)\mapsto (q,\abs{q}^2 - \omega)\,\in\,\Sigma$ is an
isomorphism between $\QH$ and $\Sigma$.

{ On the group $\QH$, the standard contact form, written as a purely
imaginary quaternion valued form, is  given by  \hspace{2mm}
$2\tilde{\Theta} = (d\omega \ - \ q \cdot d\bar q \ + \ dq\,
\cdot\bar q),$ \hspace{2mm} where $\cdot$ denotes the quaternion
multiplication,
\begin{equation}
\begin{aligned}\label{e:Heisenbegr ctct forms}
\tilde\Theta_1\ =\ \frac 12\ dx\ -\ x^{a} d t^{a}\  +\ t^{a} d
x^{a}\ -\ z^{a} d y^{a} \ +\ y^{a} d z^{a}\\
 \tilde\Theta_2\ =\ \frac 12\ dy\ -\ y^{a} d t^{a}\
+\ z^{a} d x^{a}\ + \ t^{a} d y^{a} \ -\ x^{a} d z^{a}\\
\tilde\Theta_2\ =\ \frac 12\ dz\ -\ z^{a} d t^{a}\  - \ y^{a} d x^{a}\ + \
x^{a} d y^{a} \ +\ t^{a} d z^{a}.
\end{aligned}
\end{equation}
 Since \hspace{2mm} $dp\ =\ q\cdot d\bar q\ +\ dq\,
\cdot\bar {q}\ -\ d\omega, $ \hspace{2mm} under the identification
of $\QH$ with $\Sigma$ we also have  \hspace{2mm} $2\tilde{\Theta}\
=\ - dp'\ +\ 2dq'\cdot\bar {q}'. $ \hspace{2mm}

The left invariant flat connection on $\QH$ coincides with the
Biquard connection of the qc manifold $(\QH,\tilde\Theta)$ and,
conversely, any qc manifold with flat Biquard connection is
locally isomorphic to $\QH$ \cite{IMV}.

The Cayley transform is the map
$\mathcal{C}:S\mapsto \Sigma$
from the sphere $S\ =\ \{\abs{q}^2+\abs{p}^2=1 \}\subset \Hn\times\mathbb{H}$ minus a
point to the Heisenberg group $ \Sigma$, with $\mathcal{C}$ defined by
\[
 (q', p')\ =\ \mathcal{C}\ \Big ((q, p)\Big), \qquad
 % \]
%where
%\[
q'\ =\ (1+p)^{-1} \ q, \qquad p'\ =\ (1+p)^{-1} \ (1-p)
\]
\noindent and with an inverse map $(q, p)\ =\ \mathcal{C}^{-1}\Big ((q', p')\Big)$ given
by
\[
q\ =\ \ 2(1+p')^{-1} \ q', \qquad  p\ =\  (1-p')\,(1+p')^{-1} .
\]
The unit sphere $S$  carries a natural qc structure $\tilde\eta =
dq\cdot \bar q\ +\ dp\cdot \bar p\ -\ q\cdot d\bar q -\ p\cdot d\bar
p$ which has zero torsion and  is 3-Sasakian up to a constant
factor. In \cite{IMV} it was noted that the Cayley transform is a
quaternionic contact conformal
 diffeomorphism between the quaternionic Heisenberg group with
its standard quaternionic contact structure $\tilde\Theta$ and
$S\setminus\{(-1,0)\}$ with the structure $\tilde\eta$
\begin{equation*}\label{e:Cayley transf of ctct form}
\lambda\ \cdot (\mathcal{C}_*\, \tilde\eta)\ \cdot \bar\lambda\ =\
\frac {8}{\abs{1+p'\, }^2}\, \tilde\Theta,
\end{equation*}
where $\lambda = \frac {1+p'}{\abs{1+p'\,}}$ is a unit quaternion.

\section{Curvature and the Bianchi identities}

Recall that an orthonormal frame\newline
\centerline{$\{e_1,e_2=I_1e_1,e_3=I_2e_1,e_4=I_3e_1 \dots,
e_{4n}=I_3e_{4n-3}, \xi_1, \xi_2, \xi_3 \}$}\newline
is a qc-normal frame (at a point) if the connection 1-forms of the Biquard
connection vanish (at that point). Lemma~4.5 in \cite{IMV} asserts that a
qc-normal frame exists at each point of a QC manifold.

In general, to verify any $Sp(n)Sp(1)$-invariant tensor identity at a point it
is sufficient to check it in a qc-normal
frame at that point. Further, we work in a qc-normal frame.

%\subsection{The first Bianchi identity.}

Let $b(A,B,C)$ denote the Bianchi projector,
\begin{equation}  \label{bian01}
b(A,B,C):=\sum_{(A,B,C)}\Bigl\{ (\nabla_AT)(B,C) +
T(T(A,B),C)\Bigr\}, \quad  A,B,C\in \Gamma(TM),
\end{equation}
where $\sum_{(A,B,C)}$ denotes the cyclic sum over the three tangent
vectors. With this notation the first Bianchi identity reads as
follows
\begin{equation}  \label{bian1}
\sum_{(A,B,C)}\Bigl\{R(A,B,C,D)\Bigr\}= g\Bigl ( b(A,B,C), D\Bigr )
\quad A,B,C,D\in \Gamma(TM).
\end{equation}
\begin{thrm}
\label{bianrrr} On a QC manifold the curvature of the Biquard connection satisfies the
equalities:
\begin{multline}\label{zamiana}
 R(X,Y,Z,V)-R(Z,V,X,Y)=2\sum_{s=1}^3\Big[\omega_s(X,Y)U(I_sZ,V)-
\omega_s(Z,V)U(I_sX,Y)\Big]\\
 -2\sum_{s=1}^3\Big[\omega_s(X,Z)T^0(\xi_s,Y,V)+
\omega_s(Y,V)T^0(\xi_s,Z,X)
-\omega_s(Y,Z)T^0(\xi_sX,V)-\omega_s(X,V)T^0(\xi_s,Z,Y) \Big].
\end{multline}
\begin{multline}\label{comp1}
 3R(X,Y,Z,V)-R(I_1X,I_1Y,Z,V)-R(I_2X,I_2Y,Z,V)-R(I_3X,I_3Y,Z,V)\\
= 2\Big[g(Y,Z)T^0(X,V)+g(X,V)T^0(Z,Y)-g(Z,X)T^0(Y,V)- g(V,Y)T^0(Z,X)\Big]
\\
-2\sum_{s=1}^3\Big[\omega_s(Y,Z)T^0(X,I_sV)+\omega_s(X,V)T^0(Z,I_sY)-\omega_s(Z,X)T^0(Y,I_sV)-
\omega_s(V,Y)T^0(Z,I_sX)\Big]
\\%\sum_{s=1}\Big[\omega_s(X,Y)\Big(8U(I_sZ,V)+4\rho_s(Z,V)\Big)-8\omega_s(Z,V)U(I_sX,Y)\Big]\notag\\
+\sum_{s=1}^3\Big[2\omega_s(X,Y)\Big(T^0(Z,I_sV)-T^0(I_sZ,V)\Big)-8\omega_s(Z,V)U(I_sX,Y)-\frac{Scal}{2n(n+2)}
\omega_s(X,Y)\omega_s(Z,V)\Big];
\end{multline}
\begin{multline}\label{vert1}
 R(\xi_i,X,Y,Z)= -(\nabla_XU)(I_iY,Z) +\omega_j(X,Y)\rho_k(I_iZ,\xi_i)-\omega_k(X,Y)\rho_j(I_iZ,\xi_i)\\
 -\frac14\Big[(\nabla_YT^0)(I_iZ,X)+(\nabla_YT^0)(Z,I_iX)\Big] +\frac14\Big[(\nabla_ZT^0)(I_iY,X)+
 (\nabla_ZT^0)(Y,I_iX)\Big]\\
-\omega_j(X,Z)\rho_k(I_iY,\xi_i)+\omega_k(X,Z)\rho_j(I_iY,\xi_i)
-\omega_j(Y,Z)\rho_k(I_iX,\xi_i)+\omega_k(Y,Z)\rho_j(I_iX,\xi_i)
\end{multline}
\begin{multline}\label{vert2}
 R(\xi_i,\xi_j,X,Y)=(\nabla_{\xi_i}U)(I_jX,Y)-(\nabla_{\xi_j}U)(I_iX,Y)\\
 -\frac14\Big[(\nabla_{\xi_i}T^0)(I_jX,Y)+(\nabla_{\xi_i}T^0)(X,I_jY)\Big]
 +\frac14\Big[(\nabla_{\xi_j}T^0)(I_iX,Y)+(\nabla_{\xi_j}T^0)(X,I_iY)\Big]\\
 -(\nabla_X\rho_k)(I_iY,\xi_i) -\frac{Scal}{8n(n+2)}T(\xi_k,X,Y)
 -T(\xi_j,X,e_a)T(\xi_i,e_a,Y)+T(\xi_j,e_a,Y)T(\xi_i,X,e_a),
\end{multline}
where the Ricci 2-forms are given by
\begin{equation}\label{vert023}
\begin{aligned}
3(2n+1)\rho_i(\xi_i,X)=\frac14(\nabla_{e_a}T^0)\Big[(e_a,X)-3(I_ie_a,I_iX)\Big]
-(\nabla_{e_a}U)(X,e_a)\hskip1truein\\
+\frac{2n+1}{16n(n+2)}X(Scal)\\
%\end{multline}
%\begin{multline}\label{vert024}
3(2n+1)\rho_i(I_kX,\xi_j)=-3(2n+1)\rho_i(I_jX,\xi_k)=-\frac{(2n+1)(2n-1)}{16n(n+2)}X(Scal)
\hspace{2.3cm}\\
+\frac14(\nabla_{e_a}T^0)\Big[(4n+1)(e_a,X)+3(I_ie_a,I_iX)\Big]
+2(n+1)(\nabla_{e_a}U)(X,e_a).
%\end{multline}
\end{aligned}
\end{equation}
\end{thrm}
\begin{proof}
The first Bianchi identity \eqref{bian1} yields \cite{Biq1}
\begin{multline}  \label{zam}
R(A,B,C,D)-R(C,D,A,B)=\frac12\ g\bigl(\, b(A,B,C),D\bigr ) +\frac12\
g\bigl (\,b(B,C,D),A\bigr )\\
-\frac12 \ g(b(A,C,D),B)-\frac12 \ \bigl(\,b(A,B,D),C\bigr), \quad
A,B,C,D\in \Gamma(TM).
\end{multline}
With the help of Proposition~\ref{torb} and equation \eqref{torha}
the identity \eqref{zamiana} follows.

Recall the following equality \cite[Lemma 3.8]{IMV}
\begin{equation}  \label{sp1}
R(X,Y,I_iZ,I_iV)=R(X,Y,Z,V)-2\rho_j(X,Y)\omega_j(Z,V)-2\rho_k(X,Y)%
\omega_k(Z,V).
\end{equation}
Taking into account \eqref{zamiana} and \eqref{sp1}, the properties
of the torsion listed
in Propositions~\ref{torb} and \ref{tcomp}, together with equations \eqref{propt}, %
\eqref{propu} and \eqref{t001} we find
\begin{multline}\label{cur111}
R(X,Y,Z,V)-R(I_iX,I_iY,Z,V)= 2\omega_j(X,Y)\rho_j(Z,V)+2\omega_k(X,Y)\rho_k(Z,V)\\
4\omega_j(X,Y)U(I_jZ,V)+4\omega_k(X,Y)U(I_kZ,V)-4\omega_j(Z,V)U(I_jX,Y)-4\omega_k(Z,V)U(I_kX,Y)
\\
+ 2\Big[g(Y,Z)T^0(\xi_i,I_iX,V)+g(X,V)T^0(\xi_i,I_iZ,Y)-g(Z,X)T^0(%
\xi_i,I_iY,V)- g(V,Y)T^0(\xi_i,I_iZ,X)\Big] \\
+2\Big[\omega_i(Y,Z)T^0(\xi_i,X,V)+\omega_i(X,V)T^0(\xi_i,Z,Y)-%
\omega_i(X,Z)T^0(\xi_i,Y,V)- \omega_i(Y,V)T^0(\xi_i,Z,X)\Big] \\
-\frac12\Big[\omega_j(Y,Z)\Big(T^0(X,I_jV)-T^0(I_iX,I_kV)\Big)+\omega_k(Y,Z)\Big(T^0(X,I_kV)+T^0(I_iX,I_jV)\Big)\Big]\\
-\frac12\Big[\omega_j(X,V)\Big(T^0(Y,I_jZ)-T^0(I_iY,I_kZ)\Big)+\omega_k(X,V)\Big(T^0(Y,I_kZ)+T^0(I_iY,I_jZ)\Big)\Big]\\
+\frac12\Big[\omega_j(X,Z)\Big(T^0(Y,I_jV)-T^0(I_iY,I_kV)\Big)+\omega_k(X,Z)\Big(T^0(Y,I_kV)+T^0(I_iY,I_jV)\Big)\Big]\\
+\frac12\Big[\omega_j(Y,V)\Big(T^0(X,I_jZ)-T^0(I_iX,I_kZ)\Big)+\omega_k(Y,V)\Big(T^0(X,I_kZ)+T^0(I_iX,I_jZ)\Big)\Big].
\end{multline}
Now, equality \eqref{comp1} follows from \eqref{cur111} and Theorem~\ref{sixtyseven}.

Invoking \eqref{bian01} and applying \eqref{torha} and Theorem~\ref{sixtyseven}
we have
\begin{multline}  \label{bver1}
b(\xi_i,X,Y,Z)=-(\nabla_XT)(\xi_i,Y,Z)+(\nabla_YT)(\xi_i,X,Z) \\
+2\omega_j(X,Y)\rho_k(I_iZ,\xi_i)- 2\omega_k(X,Y)\rho_j(I_iZ,\xi_i).
\end{multline}
A substitution of \eqref{bver1} in \eqref{zam} implies \eqref{vert1}.

If we take the trace in \eqref{vert1} and apply the sixth formula in Theorem~\ref{sixtyseven}
we come to
\begin{equation}\label{vert011}
n\rho_i(\xi_i,X)=\frac18(\nabla_{e_a}T^0)\Big[(e_a,X)-(I_ie_a,I_iX)\big]-\frac12\rho_k(I_jX,\xi_i)-
\frac12\rho_j(I_iX,\xi_k).
\end{equation}
Summing \eqref{vert011} and the last formula in Theorem~\ref{sixtyseven}, we obtain
\begin{equation}\label{vert012}
(n+1)\rho_i(\xi_i,X)+\frac12\rho_i(I_kX,\xi_j)=\frac18(\nabla_{e_a}T^0)\Big[(e_a,X)-(I_ie_a,I_iX)\big]
+\frac{X(Scal)}{32n(n+2)}.
\end{equation}
The second Bianchi identity
\begin{equation}  \label{secb}
\sum_{(A,B,C)}\Big\{(\nabla_AR)(B,C,D,E)+R(T(A,B),C,D,E)\Big\}=0
\end{equation}
and \eqref{torha} give
\begin{equation}  \label{bi20}
\sum_{(X,Y,Z)} \Bigl[(\nabla_XR)(Y,Z,V,W)+2\sum_{s=1}^3\omega_s(X,Y)R(%
\xi_s,Z,V,W)\Bigr]=0.
\end{equation}
We obtain  from \eqref{bi20} and \eqref{torha} that
%First, substitute $X=e_a,W=e_a$ into \eqref{bi20} and sum over $a=1,\dots,4n$. We get
\begin{multline}  \label{bi2ric}
(\nabla_{e_a}R)(X,Y,Z,e_a)-(\nabla_XRic)(Y,Z)+(\nabla_YRic)(X,Z) \\
-2\sum_{s=1}^3\Bigl[R(\xi_s,Y,Z,I_sX)-R(\xi_s,X,Z,I_sY)+\omega_s(X,Y)Ric(%
\xi_s,Z)=0.
\end{multline}
Letting $X=e_a, Y=I_ie_a$ in \eqref{bi20}  %and summing over $a=1,\dots,4n$
we find
\begin{multline}  \label{bi21}
(\nabla_{e_a}R)(I_ie_a,Z,V,W)+2n(\nabla_Z\tau_i)(V,W) \\
+2(2n-1)R(\xi_i,Z,V,W)+2R(\xi_j,I_kZ,V,W)-2R(\xi_k,I_jZ,V,W)=0.
\end{multline}
After taking the trace in \eqref{bi21} and applying the formulas in Theorem~\ref{sixtyseven}
we come to
\begin{multline}\label{vert022}
(2n-1)\rho_i(\xi_i,X)-2\rho_i(I_kX,\xi_j)=\\-\frac14\Big[(\nabla_{e_a}T^0)(e_a,X)+(\nabla_{e_a}T^0)(I_ie_a,I_iX)\Big]
-(\nabla_{e_a}U)(X,e_a)+\frac{2n-1}{16n(n+2)}X(Scal).
\end{multline}
Now, \eqref{vert012} and \eqref{vert022} yield \eqref{vert023}.% and \eqref{vert024}.

Finally, from \eqref{bian01} and an application of \eqref{torha} and Theorem~\ref{sixtyseven}
we verify that \eqref{vert2} holds.
\end{proof}

As consequence of Theorem~\ref{bianrrr} we obtain the next important Proposition.
\begin{prop}
\label{hflat} A QC manifold is locally isomorphic to the
quaternionic Heisenberg group exactly when the curvature of the
Biquard connection restricted to $H$ vanishes, $R_{|_H}=0$.
\end{prop}
\begin{proof}
Taking into account \cite[Proposition~4.11]{IMV}, in order to see
the claim it is sufficient to show that the (full) curvature
tensor vanishes. From $R_{|_H}=0$ we can conclude, cf.
\cite[Proposition~4.2, Proposition~4.3, Theorem~4.9]{IMV}, that
the vertical distribution $V$ is involutive  and
\begin{equation}  \label{van}
\rho_{t_{|_H}}=\tau_{t_{|_H}}=\zeta_{t_{|_H}}=\rho_t(\xi,.)_{|_H}=\zeta_t(%
\xi,.)_{|_H}= \tau_t(\xi,.)_{|_H} = Ric(\xi,.)_{|_H}=T(\xi_s,.)=0.
\end{equation}

Applying \eqref{van} to \eqref{vert1} and \eqref{vert2} allows us to conclude
$
R(\xi,X,Y,Z)=R(\xi_i,\xi_j,X,Y)=0$. Furthermore, \eqref{sp1curv} yields
$R(X,Y,\xi_i,\xi_j)=2\rho_k(X,Y)=0,\quad
R(X,\xi,\xi_i,\xi_j)=2\rho_k(X,\xi)=0$, and $ 4nR(\xi_s,\xi_t,\xi_i,\xi_j)=8n\rho_k(\xi_s,\xi_t)=2R(\xi_s,
\xi_t,e_a,I_ke_a)=0$, which ends the proof.
\end{proof}

\section{Quaternionic contact conformal curvature. Proof of Theorem~\ref{main1}}
In this section we define the quaternionic contact conformal curvature and prove Theorem~\ref{main1}.
\subsection{Conformal transformations}
\label{s:conf transf}

A conformal quaternionic contact transformation between two
quaternionic contact manifold is a diffeomorphism $\Phi$ which
satisfies $\Phi^*\eta=\mu\ \Psi\cdot\eta$ for some positive smooth
function $\mu$ and some matrix $\Psi\in SO(3)$ with smooth functions
as entries, where $\eta=(\eta_1,\eta_2,\eta_3)^t$ is considered as
an element of $\mathbb{R}^3$. The Biquard connection does not
change under rotations, i.e., the Biquard connection of $\Psi\cdot\eta$ and $%
\eta$ coincides. Hence, studying qc conformal transformations we may consider
only transformations $\Phi^*\eta\ =\ \mu\ \eta$.

We recall the formulas for the conformal change of the corresponding Biquard
connections from \cite{IMV}. Let $h$ be a positive smooth function on a QC
manifold $(M, \eta)$. Let $\bar\eta=\frac{1}{2h}\eta$ be a conformal
deformation of the QC structure $\eta$. We will denote the objects related
to $\bar\eta$ by over-lining the same object corresponding to $\eta$. Thus,
$d\bar\eta=-\frac{1}{2h^2}\,dh\wedge\eta\ +\ \frac{1}{2h\,}d\eta$, $\bar
g=\frac{1}{2h}g$.

The new triple $\{\bar\xi_1,\bar\xi_2,\bar\xi_3\}$, determined by
the conditions \eqref{bi1} defining the Reeb vector fields, is
$\bar\xi_s\ =\ 2h\,\xi_s\ +\ I_s\nabla h$. % where $\nabla h$ is the
%horizontal gradient defined by $g(\nabla h,X)=dh(X).$
The horizontal sub-Laplacian and the norm of the horizontal gradient are
defined respectively by $\triangle h\ =\ tr^g_H(\nabla dh)\ = \
\nabla dh(e_\alpha,e_\alpha )$, $|\nabla h|^2\ =\
dh(e_\alpha)\,dh(e_\alpha).$  The Biquard connections $\nabla$ and
$\bar\nabla$ are connected by a (1,2) tensor S,
\begin{equation}  \label{qcw2}
\bar\nabla_AB=\nabla_AB+S_AB, \qquad A,B\in\Gamma(TM).
\end{equation}
Condition \eqref{torha} yields
$
g(S_XY,Z)-g(S_YX,Z)=-h^{-1}\sum_{s=1}^3\omega_s(X,Y)dh(I_sZ),
$
while $\bar\nabla\bar g=0$ implies
$g(S_XY,Z)+g(S_XZ,Y)=-h^{-1}dh(X)g(Y,Z)$.
The last two equations determine $g(S_XY,Z)$,  %due to the equality
\begin{multline}  \label{Ivan4}
g(S_XY,Z)=-(2h)^{-1}\{dh(X)g(Y,Z)-\sum_{s=1}^3dh(I_sX)\omega_s(Y,Z)\\+
dh(Y)g(Z,X)+\sum_{s=1}^3dh(I_sY)\omega_s(Z,X)-dh(Z)g(X,Y)+%
\sum_{s=1}^3dh(I_sZ)\omega_s(X,Y)\}.
\end{multline}
Using Theorem~\ref{biqcon} we obtain after some calculations
\begin{multline}  \label{New20}
g(\bar T_{\bar\xi_1}X,Y)-2hg(T_{\xi_1}X,Y)-g(S_{\bar\xi_1}X,Y)\\
=-\nabla dh(X,I_1Y)+h^{-1}(dh(I_3X)dh(I_2Y)-dh(I_2X)dh(I_3Y)).
\end{multline}
The identity $d^2=0$ yields %\begin{equation}\label{d2}
$\nabla dh(X,Y)-\nabla dh(Y,X)=-dh(T(X,Y)).$ %\end{equation}
Applying \eqref{torha}, we have
\begin{equation}  \label{symdh}
\nabla dh(X,Y)=[\nabla dh]_{[sym]}(X,Y)-\sum_{s=1}^3 dh(\xi_s)\omega_s(X,Y),
\end{equation}
where $[.]_{[sym]}$ denotes the symmetric part of the corresponding
(0,2)-tensor.  Decomposing \eqref{New20} into [3] and [-1] parts
according to \eqref{New21}, using the properties of the torsion
tensor $T_{\xi_i}$ and \eqref{tor} we come to the next
transformation formula \cite{IMV}
\begin{gather}\label{qcw3}
g(S_{\bar\xi_i}X,Y) \ =\ -\frac{1}{4}\Big[-\nabla dh(X,I_iY)+\nabla
dh(I_iX,Y) -\nabla dh(I_jX,I_kY)+\nabla dh(I_kX,I_jY)\Bigr] \\\nonumber
\hskip.8truein -\ (2h)^{-1}\Bigl[dh(I_kX)dh(I_jY)-dh(I_jX)dh(I_kY)+
dh(I_iX)dh(Y)-dh(X)dh(I_iY)\Bigr] \\\nonumber
\hskip.11truein \ + \ \frac{1}{4n}\left(-\triangle h+2h^{-1}|\nabla
h|^2\right)\omega_i(X,Y) -dh(\xi_k)\omega_j(X,Y)+dh(\xi_j)\omega_k(X,Y).
\end{gather}

\subsection{Quaternionic contact conformal curvature}

Let $(M,g,\mathbb{Q})$ be a (4n+3)-dimensional QC manifold. We consider the
symmetric (0,2) tensor $L$ defined on $H$ by the equality
\begin{multline}  \label{lll}
L(X,Y)=\Bigl(\frac{1}{4(n+1)}Ric_{[-1]}+\frac{1}{2(2n+5)}Ric_{[3][0]}+ \frac{%
1}{32n(n+2)}Scal\,g\Bigr)(X,Y) \\
=\frac12T^0(X,Y) +U(X,Y)+\frac{Scal}{32n(n+2)}\,g(X,Y),
\end{multline}
where $Ric_{[-1]}$ is the [-1]-part of the Ricci tensor and $Ric_{[3][0]}$ is the trace-free [3]-part of $Ric$
and we use the identities in Theorem~\ref{sixtyseven} to obtain the second equality.

Let us denote the trace-free part of $L$ with $L_0$, hence
\begin{equation}  \label{l0}
L_0=\frac{1}{4(n+1)}Ric_{[-1]}+\frac{1}{2(2n+5)}Ric_{[3][0]}=\frac12T^0+U,
\end{equation}

We employ the  notation for the Kulkarni-Nomizu product of two (not
necessarily symmetric)  tensors, for example,
\begin{multline*}
(\omega_s\owedge L)(X,Y,Z,V):=\omega_s(X,Z)L(Y,V)+
\omega_s(Y,V)L(X,Z)-\omega_s(Y,Z)L(X,V)-\omega_s(X,V)L(Y,Z).
\end{multline*}
 We also note explicitly that following usual conventions we have
\[
I_s L\, (X,Y) = g(I_s L X,Y) = -L (X,I_s Y).
\]
Now, define the (0,4) tensor $WR$ on $H$ as follows
\begin{multline}\label{qcwdef}
WR(X,Y,Z,V)=R(X,Y,Z,V)+(g\owedge L)(X,Y,Z,V)+\sum_{s=1}^3(\omega _{s}\owedge
I_{s}L)(X,Y,Z,V)\\
-\frac{1}{2}\sum_{(i,j,k)}\omega_i(X,Y)\Bigl[L(Z,I_iV)-L(I_iZ,V)+L(I_jZ,I_kV)-L(I_kZ,I_jV) %
\Bigr] \\
-\sum_{s=1}^3\omega_s(Z,V)\Bigl[L(X,I_sY)-L(I_sX,Y)\Bigr]
+\frac{1}{2n}(tr L)\sum_{s=1}^3\omega_s(X,Y)\omega_s(Z,V),
\end{multline}
where $\sum_{(i,j,k)}$ denotes the cyclic sum.

A substitution of \eqref{lll} and \eqref{l0} in \eqref{qcwdef},
invoking also \eqref{propt} and \eqref{propu}, gives
\begin{multline}  \label{qcwdef1}
WR(X,Y,Z,V)= R(X,Y,Z,V)+  (g\owedge L_0)(X,Y,Z,V)+\sum_{s=1}^3(\omega_s\owedge
I_sL_0)(X,Y,Z,V)\\
-\frac12\sum_{s=1}^3\Bigl[\omega_s(X,Y)\Bigl\{T^0(Z,I_sV)-T^0(I_sZ,V)\Bigr\} %
+ \omega_s(Z,V)\Bigl\{T^0(X,I_sY)-T^0(I_sX,Y)-4U(X,I_sY)\Bigr\}\Bigr]\\
+\frac{Scal}{32n(n+2)}\Big[(g\owedge
g)(X,Y,Z,V)+\sum_{s=1}^3\Bigl((\omega_s\owedge\omega_s)(X,Y,Z,V)
+4\omega_s(X,Y)\omega_s(Z,V)\Bigr) \Big].
\end{multline}
\begin{prop}
\label{trfree} The tensor $WR$ is completely trace-free, i.e.
\begin{equation*}
Ric(WR)=\rho_s(WR)=\tau_s(WR)=\zeta_s(WR)=0.
\end{equation*}
\end{prop}
\begin{proof}
Proposition~\ref{tcomp}, \eqref{propt}, \eqref{propu} and
\eqref{lll} imply  the following identities
\begin{align}\label{t01}
T^0(\xi_s,I_sX,Y) =\frac12\Big[L(X,Y)-L(I_sX,I_sY)\Big]\hspace{6.5cm}
\\\label{u01}
U(X,Y)  =\frac14\Big[L(X,Y)+L(I_1X,I_1Y)
+L(I_2X,I_2Y)+L(I_3X,I_3Y)-\frac1n tr\,L\,g(X,Y)\Big]
%\end{align}
%\begin{multline}
\\\label{t1}
T(\xi_{i},X,Y)=
-\frac{1}{2}\Big[L(I_{i}X,Y)+L(X,I_{i}Y)\Big]+U(I_{i}X,Y)
\hspace{4.5cm}\\\notag
 =-\frac{1}{4}L(I_{i}X,Y)-\frac{3}{4}L(X,I_{i}Y)-\frac{1}{4}
L(I_{k}X,I_{j}Y)+\frac{1}{4}L(I_{j}X,I_{k}Y)-\frac{1}{4n}(tr\,L)\,g(I_{i}X,Y).
%\end{multline}
\end{align}
After a substitution of \eqref{t01} and \eqref{u01} in the first
four equations of Theorem~\ref{sixtyseven} we derive
\begin{equation}\label{ricis}
\begin{aligned}
Ric(X,Y)=\frac{2n+3}{2n}tr\,L\,g(X,Y)\hspace{7.5cm}\\+\frac{8n+11}{2}L(X,Y)+ \frac32\Bigl[
L(I_iX,I_iY)+L(I_jX,I_jY)+L(I_kX,I_kY)\Bigr]\\
%\notag \\\label{ricis}
\rho_i(X,Y)=L(X,I_iY)-L(I_iX,Y)-\frac1{2n}tr L\,\omega_i(X,Y) \hspace{4cm}\\
\tau_i(X,Y)=- \frac1n
tr\,L\,\omega_i(X,Y)\hspace{8cm}\\-\frac{n+2}{2n}\Bigl[%
L(I_iX,Y)-L(X,I_iY)+L(I_kX,I_jY)-L(I_jX,I_kY)\Bigr] \\ %\notag \\
\zeta_i(X,Y)=\frac{2n-1}{8n^2}tr\,L\,\omega_i(X,Y)\hspace{7.5cm}
\\+\frac3{8n}L(I_iX,Y)-\frac{8n+3}{8n}L(X,I_iY)+ \frac1{8n}\Bigl[%
L(I_kX,I_jY)-L(I_jX,I_kY)\Bigr].  %\notag
\end{aligned}
\end{equation}
Taking the corresponding traces in \eqref{qcwdef}, using also
\eqref{ricis}, we can  verify the claim.
\end{proof}
{ Recalling the definitions \eqref{New21} of the two $Sp(n)Sp(1)$-invariant parts of an endomorphism of $H$ and
comparing  \eqref{qcwdef1} with \eqref{comp1} we  obtain the next
Proposition.}
\begin{prop}\label{main0} On a QC manifold
the [-1]-part with respect to the first two arguments of the tensor $WR$ vanishes identically,
$$WR_{[-1]}(X,Y,Z,V)=\frac14\Big[3WR(X,Y,Z,V)-\sum_{s=1}^3WR(I_sX,I_sY,Z,V)\Big]=0.$$
{ The [3]-part with respect to the first two arguments of the tensor
$WR$ coincides with $WR$ and has the expression} 
\begin{multline}\label{qccm}
WR(X,Y,Z,V)=WR_{[3]}(X,Y,Z,V)=\frac14\Big[WR(X,Y,Z,V)+\sum_{s=1}^3WR(I_sX,I_sY,Z,V)\Big]\\
=\frac14\Big[R(X,Y,Z,V)+\sum_{s=1}^3R(I_sX,I_sY,Z,V)\Big]-\frac12\sum_{s=1}^3\omega_s(Z,V)
\Big[T^0(X,I_sY)-T^0(I_sX,Y)\Bigr]\\
+\frac{Scal}{32n(n+2)}\Big [ (g\owedge g)(X,Y,Z,V) +\sum_{s=1}^3
(\omega_s \owedge \omega_s)(X,Y,Z,V) \Big ]\\
 +(g\owedge U) (X,Y,Z,V) + \sum_{s=1}^3(\omega_s\owedge I_sU)(X,Y,Z,V).
\end{multline}
\end{prop}
\begin{dfn}
{ We denote with $W^{qc}$ the tensor $WR$ considered as a tensor of type (1,3)
with respect to the horizontal metric on $H$,
$g(W^{qc}(X,Y)Z,V)=WR(X,Y,Z,V)$
%We denote the [3]-part of the tensor $WR$ described in \eqref{qccm} by $W^{qc}, W^{qc}:=WR_{[3]}$
and call it {\textsl{the quaternionic contact conformal
curvature}}.}
\end{dfn}
\subsection{Proof of Theorem~\ref{main1}}
The relevance of $WR$ is partially justified by the following
Theorem.
\begin{thrm}
\label{qcinv} The tensor $WR$  is
{ covariant while the tensor $W^{qc}$ is invariant} under qc conformal  transformations, i.e. if
\begin{equation*}
\bar\eta=(2h)^{-1}\Psi\eta\quad {\text then} \qquad
2hWR_{\bar\eta}=WR_{\eta}, \qquad W^{qc}_{\bar\eta}=W^{qc}_{\eta}, %2h inserted after submission
\end{equation*}
for any smooth positive function $h$ and any $SO(3)$-matrix $\Psi$.
\end{thrm}
\begin{proof}
After a { standard} computation based on
\eqref{qcw2}, \eqref{Ivan4}, \eqref{qcw3}, and a suitable computer
program the relation between the curvature tensors $\bar R$ and $R$
was computed by I.Minchev \cite{IM} and presented to us in 10 pages.
After a careful study of the structure of the equation we put the
output in the following form
\begin{multline}  \label{qcw4}
2hg(\bar R(X,Y)Z,V)-g(R(X,Y)Z,V) \\
=-g\owedge M(X,Y,Z,V)-\sum_{s=1}^3 \omega_s\owedge (I_s M)(X,Y,Z,V) \\
+\frac12\sum_{(i,j,k)} \omega_i(X,Y)\Bigl[M(Z,I_iV)-M(I_iZ,V)+M(I_jZ,I_kV)-M(I_kZ,I_jV) %
\Bigr] \\
-g(Z,V)\Bigl[M(X,Y)-M(Y,X)\Bigr]+\sum_{s=1}^3\omega_s(Z,V)\Bigl[M(X,I_sY)-M(Y,I_sX)\Bigr] \\
-\frac{1}{2n}(tr M)\, \sum_{s=1}^3\omega_s(X,Y)\omega_s(Z,V) +
\frac{1}{2n}\sum_{(i,j,k)}M_i\Bigl[\omega_j(X,Y)\omega_k(Z,V) -\omega_k(X,Y)\omega_j(Z,V)%
\Bigr],
\end{multline}
where the (0,2) tensor $M$ is given by
\begin{multline}  \label{qcw5}
M(X,Y)=\frac1{2h}\Bigl(\nabla dh(X,Y)-\frac1{2h}\Bigl[dh(X)dh(Y)+
\sum_{s=1}^3dh(I_sX)dh(I_sY)+\frac12g(X,Y)|dh|^2\Bigr]\Bigr)
\end{multline}
and  $tr M=M(e_a,e_a),M_s=M(e_a,I_se_a)$
are its traces. Using \eqref{qcw5} and \eqref{symdh}, we obtain
\begin{equation}  \label{qcw6}
tr M=(2h)^{-1}\Bigl(\triangle h-(n+2)h^{-1}|dh|^2\Bigr), \qquad%\\
%M_s%(4h)^{-1}\Bigl(\nabla dh(e_a,I_se_a)-\nabla dh(I_se_a,e_a)\Bigr)%
M_s=-2n\,h^{-1}dh(\xi_s).
\end{equation}
After taking the traces in \eqref{qcw4},  using \eqref{qcw5} and the
fact that the [3]-component $(\nabla dh)_{[3]}$ of $\nabla dh$ on
$H$ is symmetric, we obtain
\begin{gather}\label{qcwric}
\overline{Ric}-Ric=4(n+1)M_{[sym]}+6M_{[3]} +\frac{2n+3}{2n}tr M\,g, \quad %\\\label{qcwscal}
\frac{\overline{Scal}}{2h}-Scal=8(n+2)tr M.
\end{gather}
The $Sp(n)Sp(1)$-invariant, [-1] and [3], parts of \eqref{qcwric}
are
\begin{gather}\label{qcwp1}
(\overline{Ric}-Ric)_{[-1]}=(n+1)M_{[sym][-1]}, \quad %\\\label{qcwp3}
(\overline{Ric}-Ric)_{[3]}=\frac{2n+5}{2}M_{[3]} + \frac{2n+3}{2n}(tr\,M)\,g.
\end{gather}
The identities in Theorem~\ref{sixtyseven}, equations \eqref{qcwric}
and \eqref{qcwp1} yield
\begin{multline}  \label{mm}
M_{[sym]}=\Bigl(\frac{1}{4(n+1)}\overline{Ric}_{[-1]}+\frac{1}{2(2n+5)}%
\overline{Ric}_{[3]}- \frac{2n+3}{32n(n+2)(2n+5)}\overline{Scal}\,\overline g%
\Bigr) \\
-\Bigl(\frac{1}{4(n+1)}Ric_{[-1]}+\frac{1}{2(2n+5)}Ric_{[3]}- \frac{2n+3}{%
32n(n+2)(2n+5)}Scal\,g\Bigr) \\
=\Bigl[\frac12\overline{T^0}+\overline{U}+\frac{\overline{Scal}}{32n(n+2)}\,%
\overline{g}\Bigr]- \Bigl[\frac12T^0+U+\frac{Scal}{32n(n+2)}\,g\Bigr].
\end{multline}
Now, from \eqref{qcw5} and \eqref{symdh} we obtain
\begin{equation}  \label{mm1}
M(X,Y)=M_{[sym]}(X,Y)-\frac{1}{2h}\sum_{s=1}^3dh(\xi_s)\omega_s(X,Y).
\end{equation}
Substituting \eqref{mm} in \eqref{mm1}, inserting the obtained equality in %
\eqref{qcw4}, and using \eqref{qcw6} completes the proof of
Theorem~\ref{qcinv}.
\end{proof}
%\begin{proof}

At this point, a combination of Theorem~\ref{qcinv} and
Proposition~\ref{main0} ends the proof of Theorem~\ref{main1} as
well.
%\end{proof}

\section{Converse problem. Proof of Theorem~\ref{main2}}

Suppose $W^{qc}=0$, hence $WR=0$. In
order to prove Theorem~\ref{main2} we search for a conformal
factor such that after a conformal transformation using this
factor the new qc structure has  Biquard connection which is flat
when restricted to the common horizontal space $H$. After we
achieve this task we can invoke Proposition~\ref{hflat} and
conclude that the given structure is locally { qc conformal to the
flat qc structure on the quaternionic Heisenberg group $\QH$.}
%qc conformally flat.
With this considerations in mind, it is then sufficient to find (locally)  a solution  $h$ of
equation \eqref{mm1} with $M_{[sym]}=-L$. In fact, a  substitution
of \eqref{mm1} in \eqref{qcw4} and an application of the condition
$W^{qc}=0=WR$ allows us to see that the qc structure
$\bar\eta=\frac1{2h}\eta$ has flat Biquard connection.

Let us consider the following overdetermined system of partial
differential equations with respect to an unknown function $u$
\begin{multline}  \label{sist1}
\nabla du(X,Y)=-du(X)du(Y)+\sum_{s=1}^3 \Bigl [
du(I_sX)du(I_sY)-du(\xi_s)\omega_s(X,Y)\Bigr ]\\
+\frac12g(X,Y)|\nabla u|^2 -L(X,Y)
\end{multline}
\begin{multline}  \label{add1}
\nabla
du(X,\xi_i)=\mathbb{B}(X,\xi_i)-L(X,I_idu)+\frac12du(I_iX)|\nabla
u|^2 \\-du(X)du(\xi_i)-du(I_jX)du(\xi_k)+du(I_kX)du(\xi_j)
\end{multline}
\begin{gather}  \label{add2}
\nabla du(\xi_i,\xi_i)=-\mathbb{B}(\xi_i,\xi_i)+\mathbb{B}%
(I_idu,\xi_i)+\frac14|\nabla u|^4-(du(\xi_i))^2+(du(\xi_j))^2+(du(\xi_k))^2, \\
\nabla du(\xi_j,\xi_i)=-\mathbb{B}(\xi_j,\xi_i)+\mathbb{B}%
(I_idu,\xi_j)-2du(\xi_i)du(\xi_j) -\frac{Scal}{16n(n+2)}du(\xi_k)
\label{add3} \\
\nabla du(\xi_k,\xi_i)=-\mathbb{B}(\xi_k,\xi_i)+\mathbb{B}(I_idu,\xi_k)
-2du(\xi_i)du(\xi_k)+\frac{Scal}{16n(n+2)}du(\xi_j).  \label{add4}
\end{gather}
Here the tensor $L$ is given by \eqref{lll}. The tensors
$\mathbb{B}(X,\xi_i)$ and $\mathbb{B}(\xi_i,\xi_j)$ do not depend on
the unknown function $u$ and { are defined in terms of $L$ and its first and second horizontal (covariant) derivatives} in \eqref{bes}
and \eqref{bst}, respectively. If we make the substitution
\begin{equation*}
2u= \ln h,\qquad 2hdu=dh, \qquad \nabla dh = 2h\nabla du +
4hdu\otimes du,
\end{equation*}
in \eqref{qcw5} we recognize that \eqref{mm1} transforms into
\eqref{sist1}. { Therefore,  our goal is to show that equation \eqref{sist1} has a solution, for which it is sufficient to verify the Ricci identities (see below). However, if equation \eqref{sist1} has a smooth solution then \eqref{add1}-\eqref{add4} appear as necessary conditions,
%However, if we consider equation \eqref{sist1} alone, we see that \eqref{add1}-\eqref{add4} appear as necessary conditions,
so we considered the complete system
\eqref{sist1}-\eqref{add4}  and reduced the question
 to showing that this system  has (locally) a smooth solution.}

The integrability conditions for the above considered
over-determined system are furnished by the Ricci identity
\begin{equation}  \label{integr}
\nabla^2du(A,B,C)-\nabla^2du(B,A,C)=-R(A,B,C,du)-\nabla du((T(A,B),C), \quad
A,B,C \in \Gamma(TM).
\end{equation}
Since \eqref{integr} is
$Sp(n)Sp(1)$-invariant it is sufficient to check it in a qc-normal frame.
%From now on
%the frame\\
%\centerline{$\{e_1,e_2=I_1e_1,e_3=I_2e_1,e_4=I_3e_1 \dots,
%e_{4n}=I_3e_{4n-3}, \xi_1, \xi_2, \xi_3 \}$}\\
%is a fixed qc-normal frame at a fixed point $p\in M$.

{ The proof of Theorem~\ref{main2} will be achieved by considering all
possible cases of \eqref{integr} and showing that the vanishing of the quaternionic
contact conformal curvature tensor $W^{qc}$ implies \eqref{integr}, which
guaranties the existence of a local smooth solution to the system
\eqref{sist1}-\eqref{add4}. The proof will be presented as a sequel
of subsections, which occupy the rest of this section. }

\subsection{Case 1, $X,Y,Z \in H$. Integrability condition ~(\ref{inte})}\hfill

When we consider equation \eqref{integr} on $H$ it takes the form
\begin{multline}  \label{inteh}
\nabla^2du(Z,X,Y)-\nabla^2du(X,Z,Y)=-R(Z,X,Y,du) \\
- 2\omega_1(Z,X)\nabla du(\xi_1,Y)-2\omega_2(Z,X)\nabla
du(\xi_2,Y)-2\omega_3(Z,X)\nabla du(\xi_3,Y),
\end{multline}
where we have used \eqref{torha}. The identity $d^2u=0$ gives
\begin{equation}  \label{comutat}
\nabla du(X,\xi_s)-\nabla du(\xi_s,X))=du(T(\xi_s,X))=T(\xi_s,X,du)
\end{equation}
After we take a covariant derivative of \eqref{sist1} along $Z\in
H$, substitute the derivatives from \eqref{sist1} and \eqref{add1},
then anti-commute the covariant derivatives, substitute the result
in \eqref{inteh} and use \eqref{qcwdef} with $WR=0$ we obtain, after
some standard calculations, that the integrability condition in this
case is
\begin{multline}  \label{inte}
(\nabla_ZL)(X,Y)-(\nabla_XL)(Z,Y)\\
=\sum_{s=1}^3\Bigl[\omega_s(Z,Y)\mathbb{B}%
(X,\xi_s)- \omega_s(X,Y)\mathbb{B}(Z,\xi_s)+2\omega_s(Z,X)\mathbb{B}%
(Y,\xi_s) \Bigr].
\end{multline}
For example, we check below that the term involving $\omega_1(Z,X)$ is $2
\mathbb{B}(Y,\xi_1)$. Indeed, the coefficient of $\omega_1(Z,X)$ in
\eqref{inteh} is calculated to be
\begin{multline*}
-\frac12\Bigl[L(Y,I_1\nabla u)-L(I_1Y,\nabla u)+L(I_2Y,I_3\nabla
u)-L(I_3Y,I_2\nabla u)\Bigr]+\frac{Scal}{16n(n+2)}\,du(I_1Y) \\
-2\nabla
du(\xi_1,Y)+du(I_1Y)|\nabla u|^2-2du(\xi_1)du(Y)-2du(\xi_3)du(I_2Y)+2du(%
\xi_2)du(I_3Y)\\
=-\frac12\Bigl[T^0(Y,I_1\nabla u)-T^0(I_1Y,\nabla u)\Bigr]+\frac{Scal}{%
16n(n+2)}\,du(I_1Y) +2L(Y,I_1\nabla
u)+2du(T(\xi_1,Y))+2\mathbb{B}(Y,\xi_1)
\\
=-\frac12\Bigl[T^0(Y,I_1\nabla u)-T^0(I_1Y,\nabla u)\Bigr] +
T^0(Y,I_1\nabla
u)+2U(Y,I_1\nabla u)+2\Bigl[du(T^0(\xi_1,Y))+U(I_1Y,\nabla u)\Bigr]\\
+2\mathbb{B}(Y,\xi_1)=2\mathbb{B}(Y,\xi_1),
\end{multline*}
where we used \eqref{comutat}, \eqref{lll} and the properties of the torsion
described in \eqref{propt},\eqref{propu} and Proposition~\ref{tcomp}.

At this point we determine the tensors $\mathbb{B}(X,\xi_s)$. Thus,
we take the traces in \eqref{inte} which give the next sequence of
equalities
\begin{equation}  \label{bes}
\begin{aligned}
&(\nabla_{e_a}L)(I_ie_a,I_iX)=(4n+1)\mathbb{B}(I_iX,\xi_i)-\mathbb{B}%
(I_jX,\xi_j)-\mathbb{B}(I_kX,\xi_k) \\
&\sum_{s=1}^3\mathbb{B}(I_sX,\xi_s)=\frac13\Bigl[(\nabla_Xtr\,L-(%
\nabla_{e_a}L)(e_a,X)\Bigr]=\frac1{4n-1}\sum_{s=1}^3(%
\nabla_{e_a}L)(I_se_a,I_sX) \\
&\mathbb{B}(X,\xi_i)=\frac1{2(2n+1)}\Bigl[(\nabla_{e_a}L)(I_ie_a,X)+\frac13%
\Bigl((\nabla_{e_a}L)(e_a,I_iX)-\nabla_{I_iX}tr\,L\Bigl)\Bigl],
\end{aligned}
\end{equation}
where the second equality in \eqref{bes} is precisely equivalent to %
\eqref{div}.

We turn to a useful technical
\begin{lemma}
\label{prost} The condition \eqref{inte} is equivalent to
\begin{equation*}
(\nabla_ZL)(X,Y)-(\nabla_XL)(Z,Y)=0\qquad {\text mod} \quad
g,\omega_1,\omega_2,\omega_3.
\end{equation*}
\end{lemma}

\begin{proof}
The condition of the lemma implies
\begin{multline}  \label{alg1}
(\nabla_ZL)(X,Y)-(\nabla_XL)(Z,Y)= g(Z,Y)C(X)-g(X,Y)C(Z)+ \\
\sum_{s=1}^3\Bigl[\omega_s(Z,Y) \mathbb{B}(X,\xi_s)- \omega_s(X,Y)\mathbb{B}%
(Z,\xi_s)+2\omega_s(Z,X) B\mathbb{(}Y,\xi_s) \Bigr],
\end{multline}
for some tensors $C(X),B(X,\xi_s)$ due to the vanishing of the
 cyclic sum $
\sum_{(Z,X,Y)}[(\nabla_ZL)(X,Y)-(\nabla_XL)(Z,Y)]=0.$ Taking traces
in \eqref{alg1} we obtain
\begin{equation*}%\label{bes1}
\begin{aligned}
&(\nabla_{e_a}L)(I_ie_a,I_iX)=(4n+1) \mathbb{B}(I_iX,\xi_i)-
\mathbb{B}
(I_jX,\xi_j)- \mathbb{B}(I_kX,\xi_k) +C(I_iX) \\
& (\nabla_{e_a}L)(e_a,X)-\nabla_{X}tr\,L=\sum_{s=1}^3(-3 \mathbb{B}
(I_sX,\xi_s)+(4n-1)C(X)) \\
& \sum_{s=1}^3(\nabla_{e_a}L)(I_se_a,I_sX)=\sum_{s=1}^3 (4n-1)
\mathbb{B}. (I_sX,\xi_s)+C(I_sX)
\end{aligned}
\end{equation*}
The last two equalities together with \eqref{div} and its consequences %
\eqref{bes} yield
\begin{equation} \label{alg2}
(4n-1)^2C(X)+3 \sum_{s=1}^3C(I_sX)=0.
\end{equation}
Solving the linear system \eqref{alg2}, we see
$((4n-1)^4+3^4)C(X)=0$. Hence, $C(X)=0$.
\end{proof}

\begin{prop}
\label{integrmain} If $W^{qc}=0$  then the
condition \eqref{inte} holds.
\end{prop}
\begin{proof}
Suppose $W^{qc}=0$, use \eqref{bi2ric} and apply \eqref{qcwdef} to
calculate
\begin{multline}\label{n=11}
(\nabla_{e_a}R)(X,Y,Z,e_a)=-(\nabla_YL)(X,Z)+(\nabla_XL)(Y,Z)\\
+\sum_{s=1}^3\Big[(\nabla_{I_sY}L)(X,I_sZ)-(\nabla_{I_sX}L)(Y,I_sZ)+(\nabla_{I_sZ}L)(X,I_sY)-
(\nabla_{I_sZ}L)(I_sX,Y)\Big]
\quad {\text mod} \quad g,\omega_s.
\end{multline}
Substituting \eqref{t01}, \eqref{u01} in \eqref{vert1} we come to
\begin{multline}\label{n=12}
-2\sum_{s=1}^3\Big[R(\xi_s,Y,Z,I_sX)-R(\xi_s,X,Z,I_sY)\Big]\\=
\sum_{s=1}^3\Big[(\nabla_{I_sY}L)(X,I_sZ)-(\nabla_{I_sX}L)(Y,I_sZ)+(\nabla_{I_sY}L)(I_sX,Z)-
(\nabla_{I_sX}L)(I_sY,Z)\Big]\\+
\frac32\sum_{s=1}^3\Big[(\nabla_{Y}L)(X,Z)-(\nabla_{X}L)(Y,Z)+(\nabla_{Y}L)(I_sX,I_sZ)-
(\nabla_{X}L)(I_sY,I_sZ)\Big] \quad {\text mod} \quad g,\omega_s.
\end{multline}
The second Bianchi identity gives
$\sum_{(X,Y,Z)}\nabla_X\rho_i(Y,Z)=0\quad {\text mod}\quad
g,\omega_s.$ Use \eqref{ricis}  to see
\begin{multline}  \label{neweq}
3\Bigl((\nabla_YL)(X,Z)-(\nabla_XL)(Y,Z)\Bigr)+\sum_{s=1}^3\Bigl(%
(\nabla_YL)(I_sX,I_sZ)-(\nabla_XL)(I_sY,I_sZ)\Bigr) \\
+ \sum_{s=1}^3\Bigl[(\nabla_{I_sZ}L)(X,I_sY)-(\nabla_{I_sZ}L)(I_sX,Y)\Bigr]%
=0 \quad {\text mod} \quad g,\omega_s.
\end{multline}
A substitution of \eqref{n=11}, \eqref{n=12}, \eqref{neweq} and
\eqref{ricis} in \eqref{bi2ric} shows, after some standard
calculations, the following identity
\begin{multline}  \label{qcbi1}
(4n+3)\Bigl[(\nabla_YL)(X,Z)-(\nabla_XL)(Y,Z)\Bigr]+
\sum_{s=1}^3\Bigl[(\nabla_{I_sY}L)(I_sX,Z)-(\nabla_{I_sX}L)(I_sY,Z)\Bigr]\\+
2\sum_{s=1}^3\Bigl[(\nabla_YL)(I_sX,I_sZ)-(\nabla_{I_sX}L)(Y,I_sZ)+\nabla_{I_sY}L)(X,I_sZ)-(\nabla_XL)(I_sY,I_sZ)
\Bigr]
=0\quad {\text mod} \quad g,\omega_s.
\end{multline}
Taking the [3]-component with respect to $X,Y$ in \eqref{qcbi1}
yields
\begin{equation}\label{n=13}
(\nabla_YL)(X,Z)-(\nabla_XL)(Y,Z)+\sum_{s=1}^3\Big[
(\nabla_{I_sY}L)(I_sX,Z)-(\nabla_{I_sX}L)(I_sY,Z)\Big]
=0\quad {\text mod} \quad g,\omega_s.
\end{equation}
A substitution of  \eqref{n=13} in \eqref{qcbi1} gives
\begin{multline}\label{qcbin1}
2n\Bigl[(\nabla_YL)(X,Z)-(\nabla_XL)(Y,Z)\Bigr]+
\sum_{s=1}^3\Bigl[(\nabla_{I_sY}L)(X,I_sZ)-(\nabla_{X}L)(I_sY,I_sZ)\Bigr]\\+
(\nabla_YL)(X,Z)-(\nabla_XL)(Y,Z)+\sum_{s=1}^3\Bigl[(\nabla_YL)(I_sX,I_sZ)-(\nabla_{I_sX}L)(Y,I_sZ)
\Bigr]
=0\quad {\text mod} \quad g,\omega_s.
\end{multline}
Taking the [-1]-component with respect to $X,Z$ of \eqref{qcbin1},
calculated with the help of \eqref{n=13}, yields
\begin{multline}  \label{qcbin2}
(6n-1)\Bigl[(\nabla_YL)(X,Z)-(\nabla_XL)(Y,Z)\Bigr]+
4\sum_{s=1}^3\Bigl[(\nabla_{I_sY}L)(X,I_sZ)-(\nabla_{X}L)(I_sY,I_sZ)\Bigr]\\
-(2n+1)\sum_{s=1}^3\Bigl[(\nabla_YL)(I_sX,I_sZ)-(\nabla_{I_sX}L)(Y,I_sZ)
\Bigr]
=0\quad {\text mod} \quad g,\omega_s.
\end{multline}
The equations \eqref{qcbin1} and \eqref{qcbin2} lead to
\begin{multline*}%\label{qcbin3}
(\nabla_YL)(X,Z)-(\nabla_XL)(Y,Z)+
\sum_{s=1}^3\Bigl[(\nabla_YL)(I_sX,I_sZ)-(\nabla_{I_sX}L)(Y,I_sZ)\Bigr]
=0\quad {\text mod} \quad g,\omega_s.
\end{multline*}
The latter  and \eqref{qcbin1} imply
\begin{equation}  \label{qcbi5}
(2n-1)\Big[(\nabla_YL)(X,Z)-(\nabla_XL)(Y,Z)\Big] =0 \quad {\text
mod} \quad g,\omega_s
\end{equation}
and Lemma~\ref{prost} completes the proof of \eqref{inte}.
\end{proof}

\subsection{Case 2, $Z,X \in H,\quad \xi_i\in V$. Integrability condition
~(\ref{inte1})}\hfill

In this case \eqref{integr} turns into
\begin{multline}  \label{intehxi}
\nabla^2du(Z,X,\xi_{i})-\nabla^2du(X,Z,\xi_{i})=-R(Z,X,\xi_{i},du)-\nabla
du(T(Z,X),\xi_{i})= \\
-2du(\xi_{j})\rho_{k}(Z,X)+2du(\xi_{k})\rho_{j}(Z,X) \\
- 2\omega_{i}(Z,X)\nabla du(\xi_{i},\xi_{i})-2\omega_{j}(Z,X)\nabla
du(\xi_{j},\xi_{i})-2\omega_{k}(Z,X)\nabla du(\xi_{k},\xi_{i}),
\end{multline}
after using \eqref{torha} and \eqref{sp1curv}.  Taking a covariant
derivative of \eqref{add1} along $Z\in H$, substituting in the
obtained equality \eqref{sist1} and \eqref{add1}, and anti-commuting
the covariant derivatives we see
\begin{multline}  \label{intver1}
\nabla^2du(Z,X,\xi_{i})-\nabla^2du(X,Z,\xi_{i})=(\nabla_Z\mathbb{B}%
)(X,\xi_{i})-(\nabla_X\mathbb{B})(Z,\xi_{i}) \\
-(\nabla_ZL)(X,I_{i}du)+(\nabla_XL)(Z,I_{i}du)
-L(X,\nabla_ZI_{i}du)+L(Z,\nabla_XI_{i}du) \\
+ {\text other \quad terms \quad coming \quad from \quad the \quad
use \quad of \quad \eqref{sist1} \quad and \quad \eqref{add1}}.
\end{multline}
Substitute \eqref{intver1} into \eqref{intehxi} use \eqref{inte} proved in Proposition~\ref{integrmain},
 also \eqref{add2}, \eqref{add3}, \eqref{add4} and the second equation in \eqref{ricis} to get after
some  calculations  that \eqref{intehxi} is equivalent to
\begin{equation}  \label{inte1}
(\nabla_Z\mathbb{B})(X,\xi_t)-(\nabla_X\mathbb{B})(Z,%
\xi_t)-L(Z,I_tL(X))+L(X,I_tL(Z))= \sum_{s=1}^32\mathbb{B}(\xi_s,\xi_t)%
\omega_s(Z,X),
\end{equation}
which is the integrability condition in this case. The
functions $\mathbb{B}(\xi_s,\xi_t)$ are uniquely determined by
%and will be determined later from the second Bianchi.
\begin{equation}  \label{bst}
\mathbb B(\xi_s,\xi_t)=\frac1{4n}\Bigl[ (\nabla_{e_a}\mathbb
B)(I_se_a,\xi_t) +L(e_a,e_b)L(I_te_a,I_se_b)\Bigr].
\end{equation}
% and do not depend on the unknown function $u$.
\begin{prop}
\label{integrmain1} If $W^{qc}=0$ then
the condition \eqref{inte1} holds.
\end{prop}
\begin{proof}
To prove the assertion it is sufficient to
 show that the left hand side of \eqref{inte1} vanishes $mod\quad \omega_s$.
Differentiating \eqref{inte} and taking the corresponding traces
yields
\begin{multline}\label{int01}
(\nabla^2_{e_a,I_{i}e_a}L)(X,Y)-(\nabla^2_{e_a,X}L)(I_{i}e_a,Y)=
-(\nabla_Y
\mathbb{B})(X,\xi_{i})-2(\nabla_X\mathbb{B})(Y,\xi_{i}) \\
+(\nabla_{I_{k}Y}\mathbb{B})(X,\xi_{j})+2(\nabla_{I_{k}X}\mathbb{B})(Y,\xi_{j})-
(\nabla_{I_{j}Y}\mathbb{B})(X,\xi_{k})-2(\nabla_{I_{j}X}\mathbb{B})(Y,\xi_{k})
\end{multline}
\begin{multline}\label{int02}
(\nabla^2_{e_a,X}L)(I_{i}e_a,Y)-(\nabla^2_{e_a,Y}L)(I_{i}e_a,X)=(\nabla_X\mathbb{
B})(Y,\xi_{i})-(\nabla_Y\mathbb{B})(X,\xi_{i}) \\
+(\nabla_{I_{k}Y}\mathbb{B})(X,\xi_{j})-(\nabla_{I_{k}X}\mathbb{B}
)(Y,\xi_{j})-(\nabla_{I_{j}Y}\mathbb{B})(X,\xi_{k})+
(\nabla_{I_{j}X}\mathbb{B} )(Y,\xi_{k})\quad {\text mod}\quad
\omega_s
\end{multline}
 \begin{multline}\label{int03}
(\nabla^2_{X,e_a}L)(I_{i}e_a,Y)=(4n+1)(\nabla_X\mathbb{B})(Y,\xi_{i})-(\nabla_X
\mathbb{B})(I_{k}Y,\xi_{j})+
(\nabla_X\mathbb{B})(I_{j}Y,\xi_{k})\quad {\text mod} \quad \omega_s
\end{multline}
\begin{multline}\label{int04}
-\nabla^2_{X,I_{i}Y}tr\,L+(\nabla^2_{X,e_a}L)(e_a,I_{i}Y)=
3(\nabla_X\mathbb{B}
)(Y,\xi_{i})-3(\nabla_X\mathbb{B})(I_{k}Y,\xi_{j})+3(\nabla_X\mathbb{B}
)(I_{j}Y,\xi_{k}).
\end{multline}
From equalities \eqref{int02} and \eqref{int03} we obtain
\begin{multline}  \label{int05}
\Bigl[\nabla^2_{X,e_a}-\nabla^2_{e_a,X}\Bigr]L(I_{i}e_a,Y)+\Bigl[%
\nabla^2_{e_a,Y}-\nabla^2_{Y,e_a}\Bigr]L(I_{i}e_a,X)\\
= 4n\Bigl[(\nabla_X\mathbb{B})(Y,\xi_{i})-(\nabla_Y\mathbb{B})(X,\xi_{i})\Bigr] \\
-\Bigl[(\nabla_X\mathbb{B})(I_{k}Y,\xi_{j})+(\nabla_{I_{k}Y}\mathbb{B}%
)(X,\xi_{j})-(\nabla_Y\mathbb{B})(I_{k}X,\xi_{j})- (\nabla_{I_{k}X}\mathbb{B}%
)(Y,\xi_{j})\Bigr] \\
+\Bigl[(\nabla_X\mathbb{B})(I_{j}Y,\xi_{k})+(\nabla_{I_{j}Y}\mathbb{B}%
)(X,\xi_{k})-(\nabla_Y\mathbb{B})(I_{j}X,\xi_{k})- (\nabla_{I_{j}X}\mathbb{B}%
)(Y,\xi_{k})\Bigr] \qquad {\text mod}\quad \omega_s.
\end{multline}
On the other hand, the Ricci identities
\begin{multline}  \label{int006}
\Bigl[\nabla^2_{X,e_a}-\nabla^2_{e_a,X}\Bigr]%
L(I_{i}e_a,Y)=-R(X,e_a,Y,e_b)L(e_b,I_{i}e_a)-4n\zeta_{i}(X,e_a)L(Y,e_a) \\
2(\nabla_{\xi_{i}}L)(X,Y)-2(\nabla_{\xi_{j}}L)(I_{k}X,Y)+2(\nabla_{%
\xi_{k}}L)(I_{j}X,Y)
\end{multline}
and the first Bianchi identity \eqref{bian1} imply
\begin{multline}  \label{int06}
\Bigl[\nabla^2_{X,e_a}-\nabla^2_{e_a,X}\Bigr]L(I_{i}e_a,Y)+\Bigl[%
\nabla^2_{e_a,Y}-\nabla^2_{Y,e_a}\Bigr]L(I_{i}e_a,X)= \\
-2\Bigl[(\nabla_{\xi_{j}}L)(I_{k}X,Y)-(\nabla_{\xi_{j}}L)(X,I_{k}Y)\Bigr] +2\Bigl[%
(\nabla_{\xi_{k}}L)(I_{j}X,Y)-(\nabla_{\xi_{k}}L)(X,I_{j}Y)\Bigr] \\
+2T(\xi_{i},Y,e_a)L(X,e_a)-2T(\xi_{j},Y,e_a)L(I_{k}X,e_a)+2T(%
\xi_{k},Y,e_a)L(I_{j}X,e_a) \\
-2T(\xi_{i},X,e_a)L(Y,e_a)+2T(\xi_{j},X,e_a)L(I_{k}Y,e_a)-2T(%
\xi_{k},X,e_a)L(I_{j}Y,e_a) \\
-R(X,Y,e_a,e_b)L(e_b,I_{i}e_a)-4n[\zeta_{i}(X,e_a)L(Y,e_a)-%
\zeta_{i}(Y,e_a)L(X,e_a)] \qquad {\text mod}\quad \omega_s.
\end{multline}
The second equality in \eqref{ricis} and a suitable contraction in
the second Bianchi identity give the next two equations valid mod $
\omega_s$
\begin{equation}  \label{xirho}
\begin{aligned}
(\nabla_{\xi_{j}}L)(X,I_{k}Y) & -(\nabla_{\xi_{j}}L)(I_{k}X,Y)=
(\nabla_{\xi_{j}}\rho_{k})(X,Y) \\
& =(\nabla_X\rho_{k})(\xi_{j},Y)-(\nabla_Y\rho_{k})(\xi_{j},X)-\rho_{k}(T(\xi_{j},X),Y)+%
\rho_{k}(T(\xi_{j},Y),X) \\
 (\nabla_{\xi_{k}}L) (X,I_{j}Y) & -(\nabla_{\xi_{k}}L)(I_{j}X,Y)=
(\nabla_{\xi_{k}}\rho_{j})(X,Y) \\
& =(\nabla_X\rho_{j})(\xi_{k},Y)-(\nabla_Y\rho_{j})(\xi_{k},X)-\rho_{j}(T(\xi_{k},X),Y)+%
\rho_{j}(T(\xi_{k},Y),X) .
\end{aligned}
\end{equation}
A substitution of \eqref{t01}, \eqref{u01} in equations
\eqref{vert023}, together with a use of \eqref{inte} and an
application of \eqref{bes} give the next
\begin{lemma}
\label{l:rho-ta} We have the following formulas for the Ricci 2-forms
\begin{equation}  \label{verf1}
\begin{aligned} \rho_{k}(\xi_{i},X) & =
\mathbb B(X,\xi_{j})-\mathbb B(I_{k}X,\xi_{i})\hskip.7truein \rho_{i}(\xi_{k},X) =
-\mathbb B(X,\xi_{j})-\mathbb B(I_{i}X,\xi_{k})\\ \rho_{i}(X,\xi_{i}) & = -\frac
{1}{4n}d(trL)(X) + \mathbb B(I_{i}X,\xi_{i}) \end{aligned}
\end{equation}
\end{lemma}

When we take the covariant derivative of \eqref{verf1}, substitute
the obtained equalities together with
 \eqref{int06}, \eqref{xirho}, %\eqref{imt07}
in \eqref{int05} we derive the formula
\begin{multline}  \label{int08}
(4n+2)\Bigl[(\nabla_X\mathbb{B})(Y,\xi_{i})-(\nabla_Y\mathbb{B})(X,\xi_{i})\Bigr]%
+ \Bigl[(\nabla_{I_{j}X}\mathbb{B})(I_{j}Y,\xi_{i})-(\nabla_{I_{j}Y}\mathbb{B}%
)(I_{j}X,\xi_{i})\Bigr] \\
+\Bigl[(\nabla_{I_{k}X}\mathbb{B})(I_{k}Y,\xi_{i})-(\nabla_{I_{k}Y}\mathbb{B}%
)(I_{k}X,\xi_{i})\Bigr]=F(X,Y) \quad {\text mod} \quad \omega_s,
\end{multline}
where the (0,2)-tensor $F$ is defined by
\begin{multline}  \label{int09}
F(X,Y)= -R(X,Y,e_a,e_b)L(e_b,I_{i}e_a)-4n\Bigl[\zeta_{i}(X,e_a)L(Y,e_a)-%
\zeta_{i}(Y,e_a)L(X,e_a)\Bigr] \\
+2T(\xi_{i},Y,e_a)L(X,e_a)-2T(\xi_{j},Y,e_a)L(I_{k}X,e_a)+2T(%
\xi_{k},Y,e_a)L(I_{j}X,e_a) \\
-2T(\xi_{i},X,e_a)L(Y,e_a)+2T(\xi_{j},X,e_a)L(I_{k}Y,e_a)-2T(%
\xi_{k},X,e_a)L(I_{j}Y,e_a) \\
+ \rho_{j}(T(\xi_{k},X),Y)-\rho_{j}(T(\xi_{k},Y),X)
+\rho_{j}(T(\xi_{k},I_{j}X),I_{j}Y)-\rho_{j}(T(\xi_{k},I_{j}Y),I_{j}X) \\
-\rho_{k}(T(\xi_{j},X),Y)+\rho_{k}(T(\xi_{j},Y),X) -
\rho_{k}(T(\xi_{j},I_{k}X),I_{k}Y)+\rho_{k}(T(\xi_{j},I_{k}Y),I_{k}X).
\end{multline}
Solving for
$(\nabla_X\mathbb{B})(Y,\xi_{i})-(\nabla_Y\mathbb{B})(X,\xi_{i})$ we
obtain
\begin{multline}  \label{fininte2}
16n(n+1)(2n+1)\Bigl[\nabla_X\mathbb{B})(Y,\xi_{i})-(\nabla_Y\mathbb{B}%
)(X,\xi_{i})\Bigr] \\
=(8n^2+8n+1)F(X,Y)+F(I_{i}X,I_{i}Y)-(2n+1)\Bigl[F(I_{j}X,I_{j}Y)+ F(I_{k}X,I_{k}Y)\Bigr] %
\quad {\text mod} \quad \omega_s.
\end{multline}
The condition $W^{qc}=0$ and \eqref{qcwdef} give
\begin{multline}  \label{curvl}
-R(X,Y,e_a,e_b)L(I_{i}e_a,e_b)=4L(X,e_a)L(Y,I_{i}e_a)-2L(X,e_a)L(I_{i}Y,e_a) \\
+2L(I_{i}X,e_a)L(Y,e_a) +2L(X,e_a)L(I_{j}Y,I_{k}e_a)-2L(I_{k}X,e_a)L(Y,I_{j}e_a) \\
-2L(X,e_a)L(I_{k}Y,I_{j}e_a)+2L(I_{j}X,e_a)L(Y,I_{k}e_a) -tr\,L\Big[%
L(X,I_{i}Y)-L(I_{i}X,Y)\Big].
\end{multline}
Using \eqref{ricis}, we get
\begin{multline}  \label{zetl}
-4n\Bigl[\zeta_{i}(X,e_a)L(Y,e_a)-\zeta_{i}(Y,e_a)L(X,e_a)\Bigr]%
=-(8n+3)L(X,e_a)L(Y,I_{i}e_a) \\
+\frac32L(X,e_a)L(I_{i}Y,e_a)-\frac32L(I_{i}X,e_a)L(Y,e_a)
-\frac12L(X,e_a)L(I_{j}Y,I_{k}e_a)+\frac12L(I_{k}X,e_a)L(Y,I_{j}e_a) \\
+\frac12L(X,e_a)L(I_{k}Y,I_{j}e_a) -\frac12L(I_{j}X,e_a)L(Y,I_{k}e_a) +\frac{2n-1}{2n%
}tr\,L\Big[L(X,I_{i}Y)-L(I_{i}X,Y)\Big].
\end{multline}
From \eqref{curvl} and \eqref{zetl} we have
\begin{multline}  \label{curzetl}
-R(X,Y,e_a,e_b)L(I_{i}e_a,e_b)-4n\Bigl[\zeta_{i}(X,e_a)L(Y,e_a)-%
\zeta_{i}(Y,e_a)L(X,e_a)\Bigr] \\
= -(8n-1)L(X,e_a)L(Y,I_{i}e_a)
-\frac12L(X,e_a)L(I_{i}Y,e_a)+\frac12L(I_{i}X,e_a)L(Y,e_a)
+\frac32L(X,e_a)L(I_{j}Y,I_{k}e_a)\\-\frac32L(I_{k}X,e_a)L(Y,I_{j}e_a)-\frac32L(X,e_a)L(I_{k}Y,I_{j}e_a)
+\frac32L(I_{j}X,e_a)L(Y,I_{k}e_a)\\
-\frac1{2n}(tr L)\Big[L(X,I_{i}Y)-L(I_{i}X,Y)\Big]
\\
=-(8n-1)L(X,e_a)L(Y,I_{i}e_a)-\frac1{2n}(tr L)\Big[L(X,I_{i}Y)-L(I_{i}X,Y)\Big] \\
+\frac12\Big[L(Y,e_a)L(I_{i}X,e_a)-L(X,e_a)L(I_{i}Y,e_a)\Big] -
\frac32\Big[L(Y,e_a)L(I_{j}X,I_{k}e_a)-L(X,e_a)L(I_{j}Y,I_{k}e_a)\Big] \\
+
\frac32\Big[L(Y,e_a)L(I_{k}X,I_{j}e_a)-L(X,e_a)L(I_{k}Y,I_{j}e_a)\Big].
\end{multline}
Since $\rho_s$ is  a (1,1)-form with respect to $I_s$, see
Proposition~\ref {sixtyseven}, we have
$$\rho_{j}(T(\xi_{k},I_{j}X),I_{j}Y)=\rho_{j}(e_a,I_{j}Y)T(\xi_{k},I_{j}X,e_a)=%
\rho_{j}(e_a,Y)T(\xi_{k},I_{j}X,I_{j}e_a).$$ Thus,  using
\eqref{ricis} we obtain the next sequence of equalities
\begin{multline}  \label{rhot}
\rho_{j}(T(\xi_{k},X),Y)+\rho_{j}(T(\xi_{k},I_{j}X),I_{j}Y)-\rho_{k}(T(\xi_{j},X),Y)-%
\rho_{k}(T(\xi_{j},I_{k}X),I_{k}Y)\\
%= \\\rho_{j}(e_a,Y)\Big[T(\xi_{k},X,e_a)+T(\xi_{k},I_{j}X,I_{j}e_a)\Big]-\rho_{k}(e_a,Y)\Big[%
%T(\xi_{j},X,e_a)+T(\xi_{j},I_{k}X,I_{k}e_a)\Big]= \\
=\Big[L(e_a,I_{j}Y)-L(I_{j}e_a,Y)-\frac1{2n}tr\,L\,\omega_{j}(e_a,Y)\Big]\Big[%
T(\xi_{k},X,e_a)+T(\xi_{k},I_{j}X,I_{j}e_a)\Big] \\
-\Big[L(e_a,I_{k}Y)-L(I_{k}e_a,Y)-\frac1{2n}tr\,L\,\omega_{k}(e_a,Y)\Big]\Big[%
T(\xi_{j},X,e_a)+T(\xi_{j},I_{k}X,I_{k}e_a)\Big]
\end{multline}
\begin{multline}  \label{rhot1}
\rho_{j}(T(\xi_{k},X),Y)+\rho_{j}(T(\xi_{k},I_{j}X),I_{j}Y)-\rho_{k}(T(\xi_{j},X),Y)-%
\rho_{k}(T(\xi_{j},I_{k}X),I_{k}Y) \\
-2T(\xi_{i},X,e_a)L(Y,e_a)+2T(\xi_{j},X,e_a)L(I_{k}Y,e_a)-2T(%
\xi_{k},X,e_a)L(I_{j}Y,e_a) \\
= L(e_a,Y)\Big[T(\xi_{k},X,I_{j}e_a)-T(\xi_{k},I_{j}X,e_a)-T(\xi_{j},X,I_{k}e_a)+T(%
\xi_{j},I_{k}X,e_a)-2T(\xi_{i},X,e_a) \Big] \\
+ L(e_a,I_{j}Y)\Big[T(\xi_{k},I_{j}X,I_{j}e_a)-T(\xi_{k},X,e_a)\Big]-L(e_a,I_{k}Y)\Big[%
T(\xi_{j},I_{k}X,I_{k}e_a)-T(\xi_{j},X,e_a)\Big] \\
-\frac1{2n}tr\,L\Big[-T(\xi_{k},X,I_{j}Y)+T(\xi_{k},I_{j}X,Y)-T(\xi_{j},I_{k}X,Y)+T(%
\xi_{j},X,I_{k}Y)\Big].
\end{multline}
The first line in \eqref{rhot1} is equal to
\begin{multline}  \label{1}
\frac1{2n}tr\,L.\,L(I_{i}X,Y)+\frac12 L(Y,e_a)\Big[5L(X,I_{i}e_a)
-L(I_{i}X,e_a)+L(I_{j}X,I_{k}e_a)-L(I_{k}X,I_{j}e_a) \Big].
\end{multline}
The second line in \eqref{rhot1} is equal to
\begin{multline}  \label{2}
\frac1{2n}tr\,L\Big[L(I_{k}X,I_{j}Y)-L(I_{j}X,I_{k}Y)\Big]+L(X,e_a)\Big[%
L(I_{k}Y,I_{j}e_a) -L(I_{j}Y,I_{k}e_a)\Big] \\
-L(I_{j}X,e_a)L(I_{j}Y,I_{i}e_a) -L(I_{k}X,e_a)L(I_{k}Y,I_{i}e_a).
\end{multline}
The third line in \eqref{rhot1} is equal to
\begin{equation}  \label{3}
-\frac1{2n}tr\,L\Big[L(I_{k}X,I_{j}Y)-L(I_{j}X,I_{k}Y)-L(X,I_{i}Y)+L(I_{i}X,Y)\Big].
\end{equation}
A substitution of \eqref{1}, \eqref{2} and \eqref{3} in
\eqref{rhot1} gives
\begin{multline}  \label{rhot11}
\rho_{j}(T(\xi_{k},X),Y)+\rho_{j}(T(\xi_{k},I_{j}X),I_{j}Y)-\rho_{k}(T(\xi_{j},X),Y)-%
\rho_{k}(T(\xi_{j},I_{k}X),I_{k}Y) \\
-2T(\xi_{i},X,e_a)L(Y,e_a)+2T(\xi_{j},X,e_a)L(I_{k}Y,e_a)-2T(%
\xi_{k},X,e_a)L(I_{j}Y,e_a) \\
= \frac1{2n}tr\,L.\,L(X,I_{i}Y)+\frac12 L(Y,e_a)\Big[5L(X,I_{i}e_a)
-L(I_{i}X,e_a)+L(I_{j}X,I_{k}e_a)-L(I_{k}X,I_{j}e_a) \Big] \\
+L(X,e_a)\Big[L(I_{k}Y,I_{j}e_a) -L(I_{j}Y,I_{k}e_a)\Big]
-L(I_{j}X,e_a)L(I_{j}Y,I_{i}e_a) -L(I_{k}X,e_a)L(I_{k}Y,I_{i}e_a).
\end{multline}
The last four lines in \eqref{int09} equal the skew symmetric sum of %
\eqref{rhot11}, which is equal to
\begin{multline}  \label{rhot111}
-5L(X,e_a)L(Y,I_{i}e_a) -\frac12\Big[L(Y,e_a)L(I_{i}X,e_a)-L(X,e_a)L(I_{i}Y,e_a)%
\Big] \\
+ \frac32\Big[L(Y,e_a)L(I_{j}X,I_{k}e_a)-L(X,e_a)L(I_{j}Y,I_{k}e_a)\Big]- \frac32%
\Big[L(Y,e_a)L(I_{k}X,I_{j}e_a)-L(X,e_a)L(I_{k}Y,I_{j}e_a)\Big] \\
+\frac1{2n}tr\,L\Big[L(X,I_{i}Y)-L(I_{i}X,Y)\Big]
-2L(I_{j}X,e_a)L(I_{j}Y,I_{i}e_a) -2L(I_{k}X,e_a)L(I_{k}Y,I_{i}e_a).
\end{multline}
A substitution of \eqref{curzetl} and \eqref{rhot111} in
\eqref{int09} yields
\begin{equation}  \label{fff}
F(X,Y)=-4(2n+1)L(X,e_a)L(Y,I_{i}e_a)-2L(I_{j}X,e_a)L(I_{j}Y,I_{i}e_a)
-2L(I_{k}X,e_a)L(I_{k}Y,I_{i}e_a).
\end{equation}
Inserting \eqref{fff} in \eqref{fininte2}  completes the proof of
\eqref{inte1}.
\end{proof}

\subsection{Case 3, $\xi\in V, \ X, Y \in H$. Integrability condition(\ref{intexih11})}\hfill

In this case \eqref{integr} reads
\begin{multline}  \label{intexih1}
\nabla^2du(\xi_{i},X,Y)-\nabla^2du(X,\xi_{i},Y)+\nabla
du(T(\xi_{i},X),Y)=-R(\xi_{i},X,Y,du).\hfill
\end{multline}
The identities  below can be used to see that the integrability
condition \eqref{intexih1} reduces to
\begin{multline}  \label{intexih11}
(\nabla_{\xi_t}L)(X,Y)+ (\nabla_X\mathbb{B})(Y,\xi_t)+L(Y,I_tL(X))+L(T(%
\xi_t,X),Y)+g(T(\xi_t,Y),L(X)) \\
= \sum_{s=1}^3\mathbb{B}(\xi_s,\xi_t)\omega_s(X,Y), \quad t=1,2,3.
\end{multline}
Notice that %that Case 3 implies Case 2, i.e.
\eqref{inte1} is  the skew-symmetric part of %a consequence of %
\eqref{intexih11}.% it is exactly the skew-symmetric part of it.

We turn to the proof of the fact the vanishing of $W^{qc}$ implies
the validity of \eqref{intexih11}. When we take a covariant
derivative along a Reeb vector field of \eqref{sist1} and a
covariant derivative along a horizontal direction of \eqref{add1}, use %
\eqref{add1}, \eqref{sist1}, \eqref{add2}, \eqref{add3}, \eqref{add4}, %
\eqref{comutat}, we see that the left hand-side of \eqref{intexih1}
equals
\begin{multline}\label{raitg}
\nabla ^{2}du(\xi _{1,}X,Y)-\nabla ^{2}du(X,\xi _{i},Y)+\nabla
du(T(\xi
_{i},X),Y) \\
=du(I_{i}Y)\left[ \mathbb{B}(I_{i}X,\xi
_{i})-\frac{1}{4n}d(tr\,L)(X)\right]
\,+du(I_{j}Y)\left[ \mathbb{B}(I_{j}X,\xi _{i})+\mathbb{B}(X,\xi _{k})%
\right] \\
+du(I_{k}Y)\left[ \mathbb{B}(I_{k}X,\xi _{i})-\mathbb{B}(X,\xi _{j})%
\right]  +g(X,Y)\mathbb{B}(du,\xi _{i})\\
-\omega _{i}(X,Y)\mathbb{B}(I_{i}du,\xi
_{i})-\omega _{j}(X,Y)\mathbb{B}(I_{j}du,\xi _{i})-\omega _{k}(X,Y)\mathbb{B}%
(I_{k}du,\xi _{i}) \\
-du(X)\mathbb{B}(Y,\xi _{i})+du(I_{i}X)\mathbb{B}(I_{i}Y,\xi _{i})+du(I_{j}X)%
\mathbb{B}(I_{j}Y,\xi _{i})+du(I_{k}X)\mathbb{B}(I_{k}Y,\xi _{i}) \\
+\frac{1}{4}(\nabla_XL)(Y,I_{i}du)-\frac{1}{4}\nabla_XL)(I_{i}Y,du)-\frac{1}{%
4}(\nabla_XL)(I_{k}Y,I_{j}du)+\frac{1}{4}(\nabla_XL)(I_{j}Y,I_{k}du) \\
-(\nabla_{\xi _{i}}L)(X,Y)\ \ -(\nabla_X\mathbb{B})(Y,\xi
_{i})+L(X,I_{i}LY)-T(\xi _{i},X,LY)-T(\xi _{i},Y,LX)\\
+\omega _{i}(X,Y)\mathbb{%
B}(\xi _{i},\xi _{i})+\omega _{j}(X,Y)\mathbb{B}(\xi _{i},\xi
_{j})+\omega _{k}(X,Y)\mathbb{B}(\xi _{i},\xi _{k}).
\end{multline}
 On the other hand,  a substitution of \eqref{t01} and \eqref{u01} in \eqref{vert1},
 and an application of \eqref{verf1}
 together with the already proven \eqref{inte} and \eqref{inte1},
 shows after  standard calculations the following equality
\begin{multline}\label{rverg}
R(\xi _{i},X,Y,Z)=\mathbb{B}(I_{j}Z,\xi _{i})\omega _{j}(X,Y)+\mathbb{B}(I_{k}Z,\xi
_{i})\omega _{k}(X,Y) \\
-\omega _{i}(Y,Z)\left[ \mathbb{B}(I_{i}X,\xi
_{i})-\frac{1}{4n}d(trL)(X) \right] -\omega _{j}(Y,Z)\Big[
\mathbb{B}(X,\xi _{k})+\mathbb{B}(I_{j}X,\xi _{i})\Big]\\ +\omega
_{k}(Y,Z)\Big[ \mathbb{B}(X,\xi _{j})-\mathbb{B} (I_{k}X,\xi
_{i})\Big]   +g(X,Z)\mathbb{B}(Y,\xi _{i})-\omega
_{i}(X,Z)\mathbb{B}(I_{i}Y,\xi _{i})\\
-\omega _{j}(X,Z)\mathbb{B}(I_{j}Y,\xi _{i})-\omega
_{k}(X,Z)\mathbb{B} (I_{k}Y,\xi _{i}) -g(X,Y)\mathbb{B}(Z,\xi
_{i})+\mathbb{B}(I_{i}Z,\xi _{i})\omega _{i}(X,Y)\\
+\frac{1}{4}\Big[ (\nabla_XL)(I_{i}Y,Z)-(\nabla_XL)(Y,I_{i}Z)+(
\nabla_XL)(I_{k}Y,I_{j}Z)-(\nabla_XL)(I_{j}Y,I_{k}Z)\Big].
\end{multline}
In the derivation of the above equation we used  the next
formulas coming  from \eqref{bes}%, and express the divergences
%of the tensor $L$  by the tensors $\mathbb{B}$.
\begin{gather}\label{e:div LI}
%\begin{aligned}
(\nabla_{e_{a,}}L)(I_{i}e_{a},X) =(4n+1)\mathbb{B}(X,\xi _{i})-\mathbb{B}%
(I_{k}X,\xi _{j})+\mathbb{B}(I_{j}X,\xi _{k}) \\
\label{e:div L} (\nabla_{e_a}L)(e_a,X) =-3\mathbb{B}(I_{i}X,\xi
_{i})-3\mathbb{B}(I_{j}X,\xi _{j})-3\mathbb{B} (I_{k}X,\xi
_{k})+d(trL)(X).
\end{gather}
Substituting  equations \eqref{rverg},  with
$Z=du$, and \eqref{raitg} in \eqref{intexih1}, we obtain \eqref{intexih11}.

%\subsection{Integrability condition \ref{intexih11}}

In the proof of the integrability condition we shall use the
following
\begin{lemma}\label{vertric-ta}
For the vertical part of the Ricci 2-forms we have the equalities
\begin{equation}\label{ric-ta-v1}
\begin{aligned}
\rho_i(\xi_j,\xi_k)& =\frac1{8n^2}(tr\,L)^2-\mathbb B(\xi_j,\xi_j)-\mathbb B(\xi_k,\xi_k)\\
\rho_i(\xi_i,\xi_j)& =\frac1{4n}d(tr\,L)(\xi_j)+\mathbb
B(\xi_i,\xi_k), \qquad \rho_i(\xi_i,\xi_k)&
=\frac1{4n}d(tr\,L)(\xi_k)-\mathbb B(\xi_i,\xi_j)
\end{aligned}
\end{equation}
\end{lemma}
\begin{proof}
From the formula for the curvature \eqref{vert2} and
Proposition~\ref{torb} it follows
\begin{align*}
 4n\rho _{i}(\xi _{i},\xi _{k}) & =(\nabla_{e_a}\rho_{j})(I_{j}e_{a},\xi _{k})+
 %\frac{trL}{n}T(\xi_{j},e_{a},I_{i}e_{a})+
 T(\xi _{i},e_{a},e_{b})T(\xi
_{k},e_{b},I_{i}e_{a})
-T(\xi_{i},e_{b},I_{i}e_{a})T(\xi _{k},e_{a},e_{b})\\
4n\rho _{j}(\xi _{i},\xi _{k})&
=-(\nabla_{e_a}\rho_{j})(I_{i}e_{a},\xi _{k})+
%\frac{trL}{n}T(\xi_{j},e_{a},I_{j}e_{a})+
T(\xi _{i},e_{a},e_{b})T(\xi _{k},e_{b},I_{j}e_{a}) -T(\xi
_{i},e_{a},I_{j}e_{b})T(\xi _{k},e_{b},e_{a}).
\end{align*}
Lemma \ref{l:rho-ta} allows us to compute the divergences
\begin{eqnarray*}
(\nabla{e_a}\rho _{i})(I_{k}e_{a},\xi _{j})&=&-
(\nabla_{a_a}\mathbb B)(I_{k}e_{a},\xi _{k})-(\nabla_{e_a}\mathbb B)(I_{j}e_{a},\xi _{j}) \\
(\nabla_{e_a}\rho _{i})(I_{i}e_{a},\xi _{j})
&=&-(\nabla_{e_a}\mathbb{B}) (e_{a},\xi
_{j})-(\nabla_{e_a}\mathbb{B})(I_{i}e_{a},\xi _{k}).
\end{eqnarray*}%
%Recalling that the torsion is completely trace free,
After a calculation in which we use the integrability condition
\eqref{inte1}, the preceding paragraphs imply the first equation of
\eqref{ric-ta-v1}. For the calculation of $\rho _{i}(\xi _{i},\xi
_{k})$ we use again \eqref{inte1} to obtain
\begin{equation*}
(\nabla_{e_a}\mathbb{B})(I_{i}e_{a},\xi _{k})
=%-L(I_{i}e_{a},I_{k}Le_{a})
-L(I_ie_a,I_ke_b)L(e_a,e_b)+4n \mathbb{B}(\xi _{i},\xi _{k}).
\end{equation*}%
Setting $s=i, Y=I_{i}X$ in \eqref{nr1}, using \eqref{ricis},
\eqref{symdh} with respect to the function $tr\,L$, together  with
Lemma~\ref{l:rho-ta} we  obtain
\begin{multline}\label{nrver1}
\Big[(\nabla_{\xi_{i}}L(X,X)+(\nabla_X\mathbb B)(X,\xi_{i})\Big] +
\Big[(\nabla_{\xi_{i}}L(I_{i}X,I_{i}X)+(\nabla_{I_{i}X}\mathbb
B)(I_{i}X,\xi_{i})\Big]=\\-\rho_{i}(e_a,X)\Big[T(\xi_{i},I_{i}X,e_a)-
T(\xi_{i},X,I_{i}e_a)\Big].
\end{multline}
Take the trace in \eqref{nrver1} and use the properties of the
torsion listed in Proposition~\ref{torb} to conclude
\begin{equation}\label{novo}
2\Big[(\nabla_{e_b}\mathbb{B})(e_{b},\xi _{i})+d(trL)(\xi
_{i})\Big]=2\rho_i(e_a,e_b)U(e_a,e_b)=0,
\end{equation}
which implies the formula for $\rho _{i}(\xi _{i},\xi _{k})$ after a
short computation.

Finally, with the help of $\rho _{i}(\xi _{i},\xi _{k})+\rho
_{j}(\xi _{j},\xi _{k})= \frac{1}{16n(n+2)}\xi
_{k}(Scal)=\frac{1}{2n}\xi _{k}(trL)$, cf. \cite[Proposition
4.4]{IMV}, we also obtain the formula for $\rho _{i}(\xi _{i},\xi
_{j})$.
\end{proof}

\begin{prop}\label{integrmain3}
If $W^{qc}=0$  then
the condition \eqref{intexih11} holds.
\end{prop}
\begin{proof}

It is sufficient to consider only the symmetric part of
\eqref{intexih11} since its skew-symmetric part is the already
established \eqref{inte1}.

Letting $A=\xi_{i}, B=X, C=Y, D=e_a, E=I_se_a$ in the second Bianchi
identity \eqref{secb} we obtain
\begin{multline}\label{nr1}
(\nabla_{\xi_{i}}\rho_s)(X,Y)-(\nabla_X\rho_s)(\xi_{i},Y)+(\nabla_Y\rho_s)(\xi_{i},X)\\
+\rho_s(T(\xi_{i},X),Y)-\rho_s(T(\xi_{i},Y),X)+2\sum_{t=1}^3\omega_t(X,Y)\rho_s(\xi_t,\xi_{i}).
\end{multline}
Setting $s=j, Y=I_{j}X$ in \eqref{nr1}, using \eqref{ricis},
Lemma~\ref{l:rho-ta}, Lemma~\ref{vertric-ta} and \eqref{inte1}, we
calculate
\begin{multline}\label{nrver2}
\Big[(\nabla_{\xi_{i}}L(X,X)+(\nabla_X\mathbb B)(X,\xi_{i})\Big] +
\Big[(\nabla_{\xi_{i}}L(I_{j}X,I_{j}X)+(\nabla_{I_{j}X}\mathbb
B)(I_{j}X,\xi_{i})\Big]=\\\rho_{j}(e_a,X)\Big[T(\xi_{i},X,I_{j}e_a)-
T(\xi_{i},I_{j}X,e_a)\Big]+\Big[(\nabla_X\mathbb
B)(I_{j}X,\xi_{k})-(\nabla_{I_{j}X}\mathbb
B)(X,\xi_{k})\Big]-2|X|^2\mathbb
B(\xi_{j},\xi_{k})\\
=2L(X,I_{k}e_a)L(I_{j}X,e_a)+\rho_{j}(e_a,X)\Big[T(\xi_{i},X,I_{j}e_a)-
T(\xi_{i},I_{j}X,e_a)\Big].
\end{multline}
Similarly, when we take $s=k, Y=I_{k}X$ in \eqref{nr1}, use
\eqref{ricis}, Lemma~\ref{l:rho-ta}, Lemma~\ref{vertric-ta} and
\eqref{inte1} it follows
\begin{multline}\label{nrver3}
\Big[(\nabla_{\xi_{i}}L(X,X)+(\nabla_X\mathbb B)(X,\xi_{i})\Big] +
\Big[(\nabla_{\xi_{i}}L(I_{k}X,I_{k}X)+(\nabla_{I_{k}X}\mathbb
B)(I_{k}X,\xi_{i})\Big]=\\\rho_{k}(e_a,X)\Big[T(\xi_{i},X,I_{k}e_a)-
T(\xi_{i},I_{k}X,e_a)\Big]-\Big[(\nabla_X\mathbb
B)(I_{k}X,\xi_{j})-(\nabla_{I_{k}X}\mathbb
B)(X,\xi_{j})\Big]+2|X|^2\mathbb
B(\xi_{j},\xi_{k})\\
=2L(I_{k}X,I_{j}e_a)L(X,e_a)+\rho_{k}(e_a,X)\Big[T(\xi_{i},X,I_{k}e_a)-
T(\xi_{i},I_{k}X,e_a)\Big].
\end{multline}
Finally, replace $X$ with $I_{i}X$ in \eqref{nrver3}, subtract the
obtained equality from \eqref{nrver2} and add the result to
\eqref{nrver1} to obtain
\begin{multline}\label{nrver4}
2\Big[(\nabla_{\xi_{i}}L(X,X)+(\nabla_X\mathbb B)(X,\xi_{i})\Big]
=2L(X,I_{k}e_a)L(I_{j}X,e_a)-2L(I_{j}X,I_{j}e_a)L(I_{i}X,e_a)\\-\rho_{i}(e_a,X)\Big[T(\xi_{i},I_{i}X,e_a)-
T(\xi_{i},X,I_{i}e_a)\Big]+\rho_{j}(e_a,X)\Big[T(\xi_{i},X,I_{j}e_a)-
T(\xi_{i},I_{j}X,e_a)\Big]\\+\rho_{k}(e_a,I_{i}X)\Big[T(\xi_{i},I_{j}X,e_a)-
T(\xi_{i},I_{i}X,I_{k}e_a)\Big].
\end{multline}
Now, using \eqref{t1} and the second equality in \eqref{ricis}
applied to \eqref{nrver4} concludes, after some standard
calculations, the  proof of  \eqref{intexih11}.
\end{proof}

\subsection{Cases 4 and 5, $\xi_{i}, \xi_{j} \in V, Y\in H$.
 Integrability conditions ~(\ref{intehxi312}), (\ref{intehxi31n}) and
(\ref{intehxi2})}\hfill

\textbf{Case 4, $\xi_{i}, \xi_{j} \in V,\quad Y\in H$}. In this case
\eqref{integr} reads
\begin{multline}  \label{intehxi2}
\nabla^2du(\xi_{i},\xi_{j},Y)-\nabla^2du(\xi_{j},\xi_{i},Y)=-R(\xi_{i},\xi_{j},Y,du)-\nabla
du(T(\xi_{i},\xi_{j}),Y).\hfill
\end{multline}
Working as in the previous case, using \eqref{add2},\eqref{add3},
\eqref{add4}, substituting \eqref{t01}, \eqref{u01} \eqref{ricis}
into \eqref{vert2}, one gets, after long standard calculations
applying the already proven \eqref{inte}, \eqref{inte1} and
\eqref{intexih11}, that \eqref{intehxi2} is equivalent to
\begin{multline}  \label{intehxi21}
(\nabla_{\xi_{i}}\mathbb{B})(X,\xi_{j})-(\nabla_{\xi_{j}}\mathbb{B}%
)(X,\xi_{i})=L(X,I_{j}e_a)\mathbb{B}(e_a,\xi_{i})-
L(X,I_{i}e_a)\mathbb{B}(e_a,\xi_{j})
\\
-L(e_a,X)\rho_{k}(I_{i}e_a,\xi_{i})-T(\xi_{i},X,e_a)\mathbb{B}(e_a,\xi_{j})+T(%
\xi_{j},X,e_a)\mathbb{B}(e_a,\xi_{i})+ \frac1n(tr L)\,\mathbb{B}(X,\xi_{k})\\
=\Big[2L(X,I_{j}e_a)+T(\xi_{j},X,e_a)\Big]\mathbb{B}(e_a,\xi_{i})-
\Big[2L(X,I_{i}e_a)+T(\xi_{i},X,e_a)\Big] \mathbb{B}(e_a,\xi_{j})+
\frac1n(tr L)\,\mathbb{B}(X,\xi_{k}).
\end{multline}
where we used Lemma~\ref{l:rho-ta} to derive the second equality.

\textbf{Case $5_a$, $X \in H,\quad \xi_{i},\xi_{j}\in V$}. In this
case \eqref{integr} becomes
\begin{multline}  \label{intehxi3}
\nabla^2du(X,\xi_{i},\xi_{j})-\nabla^2du(\xi_{i},X,\xi_{j})
=-R(X,\xi_{i},\xi_{j},du)+\nabla du(T(\xi_{i},X),\xi_{j})= \\
2du(\xi_{i})\rho_{k}(X,\xi_{i})-2du(\xi_{k})\rho_{i}(X,\xi_{i})+T(\xi_{i},X,e_a)\nabla
du(e_a,\xi_{j}).
\end{multline}
With a similar calculations as in the previous cases, we see that
\eqref{intehxi3} is equivalent to
\begin{multline}  \label{intehxi31n}
(\nabla_{\xi_{i}}\mathbb{B})(X,\xi_{j})+(\nabla_X\mathbb{B})(\xi_{i},\xi_{j}) \\
-2L(X,I_{j}e_a)\mathbb{B}(e_a,\xi_{i})+T(\xi_{i},X,e_a)\mathbb{B}%
(e_a,\xi_{j})-\frac1{2n} tr L\,\mathbb{B}(X,\xi_{k})=0.
\end{multline}
%-----------------------------------------------------------------------------------------------------------

\textbf{Case $5_b$, $X \in H,\quad \xi_{j},\xi_{j}\in V$}. In this
case \eqref{integr} reads
\begin{multline}  \label{intehxi32}
\nabla^2du(X,\xi_{j},\xi_{j})-\nabla^2
du(\xi_{j},X,\xi_{j})=-R(X,\xi_{j},\xi_{j},du)+\nabla du(T(\xi_{j},X),\xi_{j})= \\
2du(\xi_{i})\rho_{k}(X,\xi_{j})-2du(\xi_{k})\rho_{i}(X,\xi_{j})+T(\xi_{j},X,e_a)\nabla
du(e_a,\xi_{j}).
\end{multline}
and \eqref{intehxi32} is equivalent to
\begin{multline}  \label{intehxi312}
(\nabla_{\xi_{j}}\mathbb{B})(X,\xi_{j})+(\nabla_XB)(\xi_{j},\xi_{j})-2\mathbb{B}(%
\mathbf{e}_a,\xi_{j})L(X,I_{j}e_a)+T(\xi_{j},X,e_a)\mathbb{B}(e_a,\xi_{j})=0.
\end{multline}
\begin{prop}
\label{integrmain451} If $W^{qc}=0$  then
the conditions \ref{intehxi312}, \ref{intehxi31n} and \ref{intehxi2} hold.
\end{prop}
\begin{proof}
 Differentiating
the already proven \eqref{inte1} and taking the corresponding traces
we get
%for $t=1,2,3$ that
\begin{equation}  \label{in321n=1}
(\nabla^2_{X,e_a}\mathbb{B})(I_ie_a,\xi_t)+2
(\nabla_XL)(e_a,e_b)L(I_ie_a,I_te_b)=4n(\nabla_X\mathbb{B})(\xi_i,\xi_t)
\end{equation}
\begin{multline}  \label{in322n=1}
(\nabla^2_{e_a,X}\mathbb{B})(I_ie_a,\xi_t)-(\nabla^2_{e_a,I_ie_a}\mathbb{B}%
)(X,\xi_t) -2(\nabla_{\mathbf{e}_b}L)(X,I_te_a)L(I_ie_b,e_a) \\
-2(\nabla_{\mathbf{e}_b}L)(I_ie_b,e_a)L(X,I_te_a)=2(\nabla_X\mathbb{B}%
)(\xi_i,\xi_t) -2(\nabla_{I_kX}\mathbb{B})(\xi_j,\xi_t)+2(\nabla_{I_jX}%
\mathbb{B})(\xi_k,\xi_t).
\end{multline}
Subtracting \eqref{in322n=1} from \eqref{in321n=1} we obtain
\begin{multline}  \label{in323n=1}
\Bigl[\nabla^2_{X,e_a}-\nabla^2_{e_a,X}\Bigr]\mathbb{B}(I_ie_a,\xi_t)+(%
\nabla^2_{e_a,I_ie_a}\mathbb{B})(X,\xi_t) +2(\nabla_{\mathbf{e}%
_b}L)(I_ie_b,e_a)L(X,I_te_a) \\
+2\Bigl[(\nabla_XL)(e_a,e_b)-(\nabla_{\mathbf{e}_b}L)(X,e_a)\Bigr]%
L(I_ie_b,I_te_a) \\
=2(2n-1)(\nabla_X\mathbb{B})(\xi_i,\xi_t) +2(\nabla_{I_kX}\mathbb{B}%
)(\xi_j,\xi_t)-2(\nabla_{I_jX}\mathbb{B})(\xi_k,\xi_t).
\end{multline}
A use of the Ricci identities and \eqref{sp1curv} shows
\begin{multline}  \label{in324n=1}
\Bigl[\nabla^2_{X,e_a}-\nabla^2_{e_a,X}\Bigr]\mathbb{B}(I_ie_a,\xi_{i}) \\
=-R(X,e_a,I_ie_a,e_b)\mathbb{B}(e_b,\xi_{i})-R(X,e_a,\xi_{i},\xi_s)\mathbb{B}%
(I_ie_a,\xi_s)-2\omega_s(X,e_a)(\nabla_{\xi_s}\mathbb{B})(I_ie_a,\xi_{i}) \\
=-4n\zeta_i(X,,e_a)\mathbb{B}(e_a,\xi_{i})-2\rho_{k}(X,e_a)\mathbb{B}%
(I_ie_a,\xi_{j})+ 2\rho_{j}(X,e_a)\mathbb{B}(I_ie_a,\xi_{k}) \\
+2(\nabla_{\xi_i}\mathbb{B})(X,\xi_{i})-2(\nabla_{\xi_j}\mathbb{B}%
)(I_kX,\xi_{i}) +2(\nabla_{\xi_k}\mathbb{B})(I_jX,\xi_{i}).
\end{multline}
\begin{multline}  \label{in325n=1}
(\nabla^2_{e_a,I_{i}e_a}\mathbb{B})(X,\xi_{i})= \\
-\frac12\Bigl[R(e_a,I_{i}e_a,X,e_b)\mathbb{B}(e_b,\xi_{i})
+R(e_a,I_{i}e_a,\xi_{i},\xi_s)\mathbb{B}(X,\xi_s)+8n(\nabla_{\xi_{i}}\mathbb{B}%
)(X,\xi_{i})\Bigr] \\
=-2n\tau_{i}(X,e_a)\mathbb{B}(e_a,\xi_{i})-4n(\nabla_{\xi_{i}}\mathbb{B}%
)(X,\xi_{i}).
\end{multline}
Next we apply the already established \eqref{inte} and use the condition $%
L(e_a,I_se_a)=0$ to get
\begin{multline}  \label{in326n=1}
\Bigl[(\nabla_XL)(e_a,e_b)-(\nabla_{\mathbf{e}_b}L)(X,e_a)\Bigr]%
L(I_{i}e_b,I_{i}e_a) \\
=-3\mathbb{B}(e_a,\xi_{i})L(X,I_{i}e_a)+ 3\mathbb{B}(e_a,\xi_{j})L(I_{k}X,I_{i}e_a)-3%
\mathbb{B}(e_a,\xi_{k})L(I_{j}X,I_{i}e_a).
\end{multline}
\begin{multline}  \label{in327n=1}
(\nabla_{\mathbf{e}_b}L)(I_{i}e_b,e_a)L(X,I_{i}e_a) \\
=(4n+1)\mathbb{B}(e_a,\xi_{i})L(X,I_{i}e_a)- \mathbb{B}(e_a,\xi_{j})L(X,I_{j}e_a)-%
\mathbb{B}(e_a,\xi_{k})L(X,I_{k}e_a).
\end{multline}
When we substitute \eqref{in327n=1}, \eqref{in326n=1}, \eqref{in325n=1}, \eqref{in324n=1} in %
\eqref{in323n=1} we obtain
\begin{multline}  \label{int328n=1}
(1-2n)\Bigl[(\nabla_{\xi_{i}}\mathbb{B})(X,\xi_{i})+(\nabla_X\mathbb{B}%
)(\xi_{i},\xi_{i})\Bigr]- \Bigl[(\nabla_{\xi_{j}}\mathbb{B})(I_{k}X,\xi_{i})+(%
\nabla_{I_{k}X}\mathbb{B})(\xi_{j},\xi_{i})\Bigr] \\
+\Bigl[(\nabla_{\xi_{k}}\mathbb{B})(I_{j}X,\xi_{i})+(\nabla_{I_{j}X}\mathbb{B}%
)(\xi_{k},\xi_{i})\Bigr] = D_{123}(X),
\end{multline}
where $D_{ijk}(X)$ is defined by
\begin{multline}  \label{dijkn=1}
D_{ijk} (X)=\Bigl[2n\zeta_{i}(X,e_a)+n\tau_{i}(X,e_a)-(4n-2)L(X,I_{i}e_a)\Bigr]\mathbb{%
B}(e_a,\xi_{i}) \\
- \Bigl[\rho_{k}(X,I_{i}e_a)+3L(I_{k}X,I_{i}e_a)-L(X,I_{j}e_a)\Bigr]\mathbb{B}%
(e_a,\xi_{j}) \\
+ \Bigl[\rho_{j}(X,I_{i}e_a)+3L(I_{j}X,I_{i}e_a)+L(X,I_{k}e_a)\Bigr]\mathbb{B}%
(e_a,\xi_{k}).
\end{multline}
We also need the next Lemma, showing the symmetry of the vertical
tensors $\mathbb{B}$.
\begin{lemma}
\label{bijn=1} The quantities $\mathbb{B}(\xi_i,\xi_j)$ are symmetric,
\begin{equation*}
\mathbb{B}(\xi_s,\xi_t)=\mathbb{B}(\xi_t,\xi_s),\quad s,t=1,2,3.
\end{equation*}
\end{lemma}
\begin{proof}
From \eqref{bst} we obtain
\begin{equation}  \label{bst0n=1}
\mathbb{B}(\xi_{i},\xi_{j})-\mathbb{B}(\xi_{j},\xi_{i})=\frac1{4n}\Big[(\nabla_{e_a}%
\mathbb{B})(I_{i}e_a,\xi_{j})-
(\nabla_{e_a}\mathbb{B})(I_{j}e_a,\xi_{i})\Big].
\end{equation}
On the other hand, \eqref{bes} imply
\begin{multline}  \label{bst1n=1}
2(2n+1)(4n-1)\Big[\mathbb{B}(I_{i}X,\xi_{j})- \mathbb{B}(I_{j}X,\xi_{i})\Big] \\
=(4n+1)\Big[(\nabla_{e_a}L)(I_{i}e_{e_a},I_{j}X)-(\nabla_{e_a}L)(I_{j}e_{e_a},I_{i}X)%
\Big]+2(\nabla_{e_a}L)(I_{k}e_a,X).
\end{multline}
Substitute \eqref{bst1n=1} into \eqref{bst0n=1} to get
\begin{multline}  \label{bst2n=1}
8n(2n+1)(4n-1)\Big[\mathbb{B}(\xi_{i},\xi_{j})-\mathbb{B}(\xi_{j},\xi_{i})\Big] \\
=(4n+1)\Big[(\nabla^2_{e_b,e_a}L)(I_{i}e_{e_a},I_{j}e_b)-(%
\nabla^2_{e_b,e_a}L)(I_{j}e_{e_a},I_{i}e_b)\Big]+
2(\nabla^2_{e_b,e_a}L)(I_{k}e_a,e_b).
\end{multline}
We calculate using \eqref{inte} and \eqref{bst0n=1} that
\begin{multline}  \label{bst3n=1}
(\nabla^2_{e_b,e_a}L)(I_{k}e_a,e_b)=(4n+1)(\nabla_{e_a}\mathbb{B}%
)(e_a,\xi_{k})+(\nabla_{e_a}\mathbb{B})(I_{i}e_a,\xi_{j}) -(\nabla_{e_a}\mathbb{B}%
)(I_{j}e_a,\xi_{i}) \\
= (4n+1)(\nabla_{e_a}\mathbb{B})(e_a,\xi_{k})+4n\Big[\mathbb{B}(\xi_{i},\xi_{j})-%
\mathbb{B}(\xi_{j},\xi_{i}) \Big].
\end{multline}
The Ricci identities, the symmetry of $L$ and \eqref{ricis} imply
\begin{multline}  \label{bst4n=1}
\Big[(\nabla^2_{e_b,e_a}L)(I_{i}e_{e_a},I_{j}e_b)-(%
\nabla^2_{e_b,e_a}L)(I_{j}e_{e_a},I_{i}e_b)\Big] \\
=
\zeta_{j}(e_b,e_a)L(e_a,I_{i}e_b)-\zeta_{i}(e_b,e_a)L(e_a,I_{j}e_b)+2%
\omega_s(e_b,e_a)(\nabla_{\xi_{k}}L)(I_{j}e_a,I_{i}e_b)
=2\nabla_{\xi_{k}}tr\,L.
\end{multline}
Substitute \eqref{bst4n=1} and \eqref{bst3n=1} in \eqref{bst2n=1}
and apply \eqref{novo} to conclude
\begin{equation*}  \label{bst5n=1}
8n(4n^2+n-1)\Big[\mathbb{B}(\xi_{i},\xi_{j})-\mathbb{B}(\xi_{j},\xi_{i})\Big]=(4n+1)%
\Big[\nabla_{\xi_{k}}tr\,L
+(\nabla_{e_a}\mathbb{B})(e_a,\xi_{k})\Big]=0.
\end{equation*}
\end{proof}

The second Bianchi identity \eqref{secb} taken with respect to
$A=\xi_{i}, B=\xi_{j},C=X, D=e_a, E=I_se_a$ and the formulas
described in Theorem~\ref{sixtyseven} yield
\begin{multline}\label{rhon=1}
(\nabla_{\xi_{i}}\rho_s)(\xi_{j},X)-(\nabla_{\xi_{j}}\rho_s)(\xi_{i},X)+(\nabla_X\rho_s)(\xi_{i},\xi_{j})\\
=\rho_s(T(\xi_{i},X),\xi_{j})-\rho_s(T(\xi_{j},X),\xi_{i})+\rho_s(e_a,X)
\rho_{k}(I_{i}e_a.\xi_{i})+\frac{tr\,L}n\rho_s(\xi_{k},X).
\end{multline}
Setting successively  $s=1,2,3$  in \eqref{rhon=1}, using
\eqref{comutat} with respect to the function $tr\,L$ and applying
Lemma~\ref{l:rho-ta} and Lemma~\ref{vertric-ta}, we obtain after
some calculations
\begin{equation}\label{n=123}
\begin{aligned}
\Big[(\nabla_{\xi_{i}}\mathbb
B)(I_{i}X,\xi_{j})-(\nabla_{\xi_{j}}\mathbb B)(I_{i}X,\xi_{i})\Big]
-\Big[(\nabla_{\xi_{i}}\mathbb B)(X,\xi_{k})+(\nabla_X\mathbb B)(\xi_{i},\xi_{k})\Big]=\alpha_{ijk}(X)\\
\Big[(\nabla_{\xi_{i}}\mathbb
B)(I_{j}X,\xi_{j})-(\nabla_{\xi_{j}}\mathbb B)(I_{j}X,\xi_{i})\Big]
-\Big[(\nabla_{\xi_{j}}\mathbb B)(X,\xi_{k})+(\nabla_X\mathbb B)(\xi_{j},\xi_{k})\Big]=\beta_{ijk}(X)\\
\Big[(\nabla_{\xi_{i}}\mathbb
B)(I_{k}X,\xi_{j})-(\nabla_{\xi_{j}}\mathbb B)(I_{k}X,\xi_{i})\Big]
+\Big[(\nabla_{\xi_{i}}\mathbb B)(X,\xi_{i})+(\nabla_X\mathbb B)(\xi_{i},\xi_{i})\Big]\qquad \qquad\quad \\
+\Big[(\nabla_{\xi_{j}}\mathbb B)(X,\xi_{j})+(\nabla_X\mathbb
B)(\xi_{j},\xi_{j})\Big]=\gamma_{ijk}(X),
\end{aligned}
\end{equation}
where
\begin{equation}\label{alpha}
\begin{aligned}
\alpha_{ijk}(X)=\rho_{i}(e_a,\xi_{i})T(\xi_{j},X,e_a)& -\rho_{i}(e_a,\xi_{j})
T(\xi_{i},X,e_a)-\rho_{i}(e_a,X)\rho_{k}(I_{i}e_a,\xi_{i})\\
& +\frac1{4n}d(tr\,L)(e_a)T(\xi_{j},X,e_a)-\frac{tr\,L}{n}\rho_{i}(\xi_{k},X)\\
\beta_{ijk}(X)=\rho_{j}(e_a,\xi_{i})T(\xi_{j},X,e_a)& -\rho_{j}(e_a,\xi_{j})
T(\xi_{i},X,e_a)-\rho_{j}(e_a,X)\rho_{k}(I_{i}e_a,\xi_{i})\\
& -\frac1{4n}d(tr\,L)(e_a)T(\xi_{i},X,e_a)-\frac{tr\,L}{n}\rho_{j}(\xi_{k},X)\\
\gamma_{ijk}(X)=\rho_{k}(e_a,\xi_{i})T(\xi_{j},X,e_a)& -\rho_{k}(e_a,\xi_{j})
T(\xi_{i},X,e_a)-\rho_{k}(e_a,X)\rho_{k}(I_{i}e_a,\xi_{i})\\
& +\frac1{4n^2}(tr\,L)d(tr\,L)(X)-\frac{tr\,L}{n}\rho_{k}(\xi_{k},X)
\end{aligned}
\end{equation}
Now we can solve the system consisting of \eqref{dijkn=1} and
\eqref{n=123}. Indeed, \eqref{dijkn=1} and Lemma~\ref{bijn=1} imply
\begin{multline}\label{nn=12}
(1-2n)\Big[[(\nabla_{\xi_{i}}\mathbb B)(X,\xi_{i})+(\nabla_X\mathbb
B)(\xi_{i},\xi_{i})]
+[(\nabla_{\xi_{j}}\mathbb B)(X,\xi_{j})+(\nabla_X\mathbb B)(\xi_{j},\xi_{j})]\Big]\\
\Big[(\nabla_{\xi_{i}}\mathbb
B)(I_{k}X,\xi_{j})-(\nabla_{\xi_{j}}\mathbb B)(I_{k}X,\xi_{i})\Big]
+\Big[(\nabla_{\xi_{k}}\mathbb B)(I_{j}X,\xi_{i})+(\nabla_{I_{j}X}\mathbb B)(\xi_{k},\xi_{i})\Big]\\
-\Big[(\nabla_{\xi_{k}}\mathbb
B)(I_{i}X,\xi_{j})+(\nabla_{I_{i}X}\mathbb
B)(\xi_{k},\xi_{j})\Big]=D_{123}(X)+D_{231}(X).
\end{multline}
The last identity in \eqref{n=123} and \eqref{nn=12} yields
\begin{multline}\label{nn=13}
2n\Big[(\nabla_{\xi_{i}}\mathbb
B)(I_{k}X,\xi_{j})-(\nabla_{\xi_{j}}\mathbb B)(I_{k}X,\xi_{i})\Big]
+\Big[(\nabla_{\xi_{k}}\mathbb B)(I_{j}X,\xi_{i})+(\nabla_{I_{j}X}\mathbb B)(\xi_{k},\xi_{i})\Big]\\
-\Big[(\nabla_{\xi_{k}}\mathbb
B)(I_{i}X,\xi_{j})+(\nabla_{I_{i}X}\mathbb
B)(\xi_{k},\xi_{j})\Big]=D_{ijk}(X)+D_{231}(X)+(2n-1)\gamma_{ijk}(X).
\end{multline}
The first two equalities in \eqref{n=123} together with \eqref{nn=13} lead to
\begin{multline}\label{nn=14}
2(n+1)\Big[(\nabla_{\xi_{i}}\mathbb
B)(I_{k}X,\xi_{j})-(\nabla_{\xi_{j}}\mathbb B)(I_{k}X,\xi_{i})\Big]
+\Big[(\nabla_{\xi_{k}}\mathbb
B)(I_{j}X,\xi_{i})-(\nabla_{\xi_{i}}\mathbb B)(I_{j}X,\xi_{k})\Big]
+\\\Big[(\nabla_{\xi_{j}}\mathbb
B)(I_{i}X,\xi_{k})-(\nabla_{\xi_{k}}\mathbb
B)(I_{i}X,\xi_{j})\Big]=A_{ijk}(X),
\end{multline}
where
\begin{equation}\label{nn=15}
A_{ijk}(X)=D_{ijk}(X)+D_{jki}(X)+(2n-1)\gamma_{ijk}(X)+\alpha_{ijk}(I_{j}X)-\beta_{ijk}(I_{i}X).
\end{equation}
Consequently, we derive easily that
\begin{multline}\label{nn=16}
2(n+2)(2n+1)\Big[(\nabla_{\xi_{i}}\mathbb
B)(I_{k}X,\xi_{j})-(\nabla_{\xi_{j}}\mathbb
B)(I_{k}X,\xi_{i})\Big]\\=(2n+3)A_{ijk}(X)- A_{jki}(X)-A_{kij}(X).
\end{multline}
The second equality in \eqref{ricis} together with \eqref{t1} and Lemma~\ref{l:rho-ta} applied to \eqref{alpha} and
\eqref{dijkn=1}, after standard calculations, give
\begin{multline}\label{nn=17}
\alpha _{ijk}(I_{j}X)-\beta
_{ijk}(I_{i}X)=\frac{1}{2}L(X,e_{a})\Big[ \mathbb B(I_{i}e_{a},\xi
_{i})+\mathbb B(I_{j}e_{a},\xi _{j})+\mathbb B(I_{k}e_{a},\xi _{k})\Big]
\\
+\frac{1}{2}L(I_{i}X,e_{a})\Big[ -\mathbb B(e_{a},\xi
_{i})-3\mathbb B(I_{k}e_{a},\xi _{j})-2\mathbb B(I_{j}e_{a},\xi
_{k})\Big] \\
+\frac{1}{2}L(I_{j}X,e_{a})\Big[ 3\mathbb B(I_{k}e_{a},\xi _{i})-\mathbb B(e_{a},\xi
_{j})+2\mathbb B(I_{i}e_{a},\xi _{k})\Big] \\
+\frac{1}{2}L(I_{k}X,e_{a})\Big[-5\mathbb B(I_{j}e_{a},\xi
_{i})+5\mathbb B(I_{i}e_{a},\xi _{j}) -\mathbb B(e_{a},\xi _{k})\Big]
+\frac{3}{2n}\left( tr\,L\right) \mathbb B(I_{k}X,\xi _{k})
\end{multline}
\begin{multline}\label{nn=18}
\gamma _{ijk}(X)=
=\frac{1}{2n}\left( tr\,L\right) \mathbb B(I_{k}X,\xi _{k})\\-\frac{5}{2}L(X,e_{a})%
\Big[ \mathbb B(I_{i}e_{a},\xi _{i})+\mathbb B(I_{j}e_{a},\xi _{j})\Big] -\frac{3}{2}%
L(I_{k}X,e_{a})\Big[ \mathbb B(I_{j}e_{a},\xi _{i})-\mathbb B(I_{i}e_{a},\xi
_{j})\Big]
\end{multline}
\begin{multline}\label{nn=19}
D_{ijk}(X)+D_{jki}(X)+(2n-1)\gamma
_{ijk}(X)=\frac{2n+1}{4}L(X,e_{a})\Big[ \mathbb B(I_{i}e_{a},\xi
_{i})+\mathbb B(I_{j}e_{a},\xi
_{j})\Big] \\
+\frac{1}{4}L(I_{i}X,e_{a})\Big[ -\left( 2n+1\right) \mathbb B(e_{a},\xi
_{i})+\left( 2n+3\right) \mathbb B(I_{k}e_{a},\xi _{j})+8\mathbb B(I_{j}e_{a},\xi _{k})%
\Big] \\
+\frac{1}{4}L(I_{j}X,e_{a})\Big[ -\left( 2n+3\right)
\mathbb B(I_{k}e_{a},\xi
_{i})-\left( 2n+1\right) \mathbb B(e_{a},\xi _{j})-8\mathbb B(I_{i}e_{a},\xi _{k})\Big] \\
+\frac{1-10n}{4}L(I_{k}X,e_{a})\Big[ \mathbb B(I_{j}e_{a},\xi
_{i})-\mathbb B(I_{i}e_{a},\xi _{j})\Big]\\
+\frac{1}{4n}\left( trL\right) %
\Big[ (1-2n)\,\mathbb{B}(I_{i}X,\xi _{i})+\left( 1-2n\,\right) \mathbb{B}%
(I_{j}X,\xi _{j})+8n\,\mathbb{B}(I_{k}X,\xi _{k})\Big].
\end{multline}
A substitution of \eqref{nn=17}, \eqref{nn=18} and \eqref{nn=19} in
\eqref{nn=15} shows
\begin{multline}\label{nn=110}
A_{ijk}(X)=\frac{1}{4n}\left( tr\,L\right) \left[
(1-2n)\,\mathbb{B}(I_{i}X,\xi
_{i})+\left( 1-2n\,\right) \mathbb{B}(I_{j}X,\xi _{j})+\left( 8n+6\right) \,%
\mathbb{B}(I_{k}X,\xi _{k})\right]  \\
+\frac{1}{4}L(X,e_{a})\left[ \left( 2n+3\right) \mathbb B(I_{i}e_{a},\xi
_{i})+\left( 2n+3\right) \mathbb B(I_{j}e_{a},\xi _{j})+2\mathbb B(I_{k}e_{a},\xi _{k})%
\right]  \\
+\frac{1}{4}L(I_{i}X,e_{a})\left[ -\left( 2n+3\right) \mathbb B(e_{a},\xi
_{i})+\left( 2n-3\right) \mathbb B(I_{k}e_{a},\xi _{j})+4\mathbb B(I_{j}e_{a},\xi _{k})%
\right]  \\
+\frac{1}{4}L(I_{j}X,e_{a})\left[ -\left( 2n-3\right)
\mathbb B(I_{k}e_{a},\xi
_{i})-\left( 2n+3\right) \mathbb B(e_{a},\xi _{j})-4\mathbb B(I_{i}e_{a},\xi _{k})\right]  \\
+\frac{1}{4}L(I_{k}X,e_{a})\left[ -\left( 10n+9\right)
\mathbb B(I_{j}e_{a},\xi _{i})+\left( 10n+9\right) \mathbb B(I_{i}e_{a},\xi
_{j})-2\mathbb B(e_{a},\xi _{k})\right]
\end{multline}
Plugging \eqref{nn=110} in \eqref{nn=16} and using \eqref{t1} we
obtain
\begin{multline}  \label{intehxi21n1}
(\nabla_{\xi_{i}}\mathbb{B})(I_{k}X,\xi_{j})-(\nabla_{\xi_{j}}\mathbb{B})(I_{k}X,\xi_{i})=\frac1n(tr
L)\,\mathbb{B}(I_{k}X,\xi_{k})
\\+\Big[2L(I_{k}X,I_{j}e_a)+T(\xi_{j},I_{k}X,e_a)\Big]\mathbb{B}(e_a,\xi_{i})- \Big[2L(I_{k}X,I_{i}e_a)+T(\xi_{i},I_{k}X,e_a)\Big]
\mathbb{B}(e_a,\xi_{j}).
\end{multline}
Hence, \eqref{intehxi21} follows. Substituting \eqref{intehxi21} in
the first equality of \eqref{n=123} we obtain \eqref{intehxi31n}.
Inserting \eqref{intehxi31n} in \eqref{int328n=1} we see
\eqref{intehxi32}.
\end{proof}

\subsection{Case 6, $\xi_{k},\xi_{i},\xi_{j}\in V$. Integrability conditions ~(\ref{intehxi41}) and
(\ref{intehxi441})}\hfill

\textbf{Case $6_{a}$, $\xi_{k},\xi_{i},\xi_{j}\in V$}. In this case
the Ricci identity  \eqref{integr} becomes
\begin{multline}  \label{intehxi4}
\nabla^2du(\xi_{k},\xi_{i},\xi_{j})-
\nabla^2du(\xi_{i},\xi_{k},\xi_{j})=-R(\xi_{k},\xi_{i},\xi_{j},du)-\nabla
du(T(\xi_{k},\xi_{i}),\xi_{j}) \\
=2du(\xi_{i})\rho_{k}(\xi_{k},\xi_{i})-2du(\xi_{k})\rho_{i}(\xi_{k},\xi_{i})+\rho_{j}(I_{k}e_a,%
\xi_{k})\nabla du(e_a,\xi_{j}) +\frac1n tr\,L\,\nabla
du(\xi_{j},\xi_{j}).
\end{multline}
After some calculations we see that \eqref{intehxi4} is equivalent
to
\begin{multline*}
(\nabla_{\xi_{i}}\mathbb{B})(\xi_{j},\xi_{k})-(\nabla_{\xi_{k}}\mathbb{B}%
)(\xi_{j},\xi_{i})+ \mathbb{B}(I_{k}e_a,\xi_{j})\mathbb{B}(e_a,\xi_{i})-\mathbb{B}%
(I_{i}e_a,\xi_{j})\mathbb{B}(e_a,\xi_{k}) \\
-\rho_{i}(I_{k}e_a,\xi_{k})\mathbb{B}(e_a,\xi_{i})+\rho_{k}(I_{i}e_a,\xi_{i})\mathbb{B}%
(e_a,\xi_{k})-\rho_{j}(I_{k}e_a,\xi_{k})\mathbb{B}(e_a,\xi_{j}) \\
+\frac1{2n}\Bigl[trL\,\mathbb{B}(\xi_{i},\xi_{i})- trL\,\mathbb{B}%
(\xi_{k},\xi_{k})-2trL\,\mathbb{B}(\xi_{j},\xi_{j})\Bigr]=0.
\end{multline*}
Using Lemma \ref{l:rho-ta} and the above equation shows that the
integrability condition in this case is
\begin{multline}\label{intehxi41}
(\nabla_{\xi _{i}}\mathbb{B})(\xi _{k},\xi _{j})-(\nabla_{\xi _{k}}
\mathbb{B})(\xi _{i},\xi _{j}) =\frac{1}{2n}(tr\,L\,)\left[
\mathbb{B}(\xi _{i},\xi _{i})-2\mathbb{B}(\xi _{j},\xi
_{j})+\mathbb{B}(\xi _{k},\xi _{k}) \right]\\
+2\mathbb{B}(e_{a},\xi _{i})\mathbb{B}(I_{j}e_{a},\xi
_{k})+\mathbb{B} (e_{a},\xi _{i})\mathbb{B}(I_{k}e_{a},\xi
_{j})+\mathbb{B}(I_{i}e_{a},\xi _{k})\mathbb{B}(e_{a},\xi _{j}).
\end{multline}

\textbf{Case $6_{b}$, $\xi_{k},\xi_{j},\xi_{j}\in V$}. Here,
equation \eqref{integr} reads
\begin{multline}  \label{intehxi44}
\nabla^2du(\xi_{k},\xi_{j},\xi_{j})-\nabla^2du(\xi_{j},\xi_{k},\xi_{j})=-R(\xi_{k},\xi_{j},\xi_{j},du)-\nabla
du(T(\xi_{k},\xi_{j}),\xi_{j})\\=
2du(\xi_{i})\rho_{k}(\xi_{k},\xi_{j})-2du(\xi_{k})\rho_{i}(\xi_{k},\xi_{j})-\rho_{i}(I_{k}e_a,%
\xi_{k})\nabla du(e_a,\xi_{j}) -\frac1n (tr\,L)\,\nabla
du(\xi_{i},\xi_{j}).
\end{multline}
 A small calculation  shows that \eqref{intehxi44} is equivalent to
\begin{multline}\label{intehxi441}
(\nabla_{\xi _{j}} \mathbb{B})(\xi _{k},\xi _{j})-(\nabla_{\xi _{k}}
\mathbb{B})(\xi _{j},\xi _{j})\\
=-\mathbb B(I_{k}e_{a},\xi _{j})\mathbb{B}(e_{a},\xi _{j})
+3\mathbb{B}(I_{j}e_{a},\xi _{k})\mathbb{B}(e_{a},\xi _{j})
+\frac{3}{2n}(tr\,L)\, \mathbb{B}(\xi _{i},\xi _{j})=0.
\end{multline}
\begin{prop}
\label{integrmain6} If $W^{qc}=0$  then
the conditions \ref{intehxi41}, \ref{intehxi441} hold.
\end{prop}
\begin{proof}
Differentiate
\eqref{intehxi31n} and take the corresponding trace to get
\begin{multline}
(\nabla _{e_{a},\xi _{i}}^{2}\mathbb{B})(I_{k}e_{a},\xi
_{j})+(\nabla
_{e_{a},I_{k}e_{a}}^{2}\mathbb{B})(\xi _{i},\xi _{j})=  \label{int60} \\
2(\nabla _{e_{b}}L)(I_{k}e_{b},I_{j}e_{a})\mathbb{B}(e_{a},\xi
_{i})+2L(I_{k}e_{b},I_{j}e_{a})(\nabla _{e_{b}}\mathbb{B})(e_{a},\xi _{i}) \\
-(\nabla _{e_{b}}T)(\xi _{i},I_{k}e_{b},e_{a})\mathbb{B}(e_{a},\xi
_{j})-T(\xi _{i},I_{k}e_{b},e_{a})(\nabla
_{e_{b}}\mathbb{B})(e_{a},\xi _{j})
\\
+\frac{1}{2n}d(tr\,L)(e_{a})\mathbb{B}(I_{k}e_{a},\xi _{k})+\frac{1}{2n}%
(tr\,L)(\nabla _{e_{a}}\mathbb{B})(I_{k}e_{a},\xi _{k}).
\end{multline}%
On the other hand, the Ricci identities, \eqref{bst}, \eqref{sp1curv} and %
\eqref{ricis} yield
\begin{equation}
(\nabla _{e_{a},I_{k}e_{a}}^{2}\mathbb{B})(\xi _{i},\xi
_{j})=-4n(\nabla _{\xi _{k}}\mathbb{B})(\xi _{i},\xi
_{j})+4(tr\,L)\mathbb{B}(\xi _{j},\xi _{j})-4(tr\,L)\mathbb{B}(\xi
_{i},\xi _{i}).  \label{int601}
\end{equation}%
\begin{multline}
(\nabla _{e_{a},\xi _{i}}^{2}\mathbb{B})(I_{k}e_{a},\xi
_{j})=(\nabla _{\xi
_{i},e_{a}}^{2}\mathbb{B})(I_{k}e_{a},\xi _{j})+4n\zeta _{k}(\xi _{i},e_{a})%
\mathbb{B}(e_{a},\xi _{j})  \label{int602} \\
-2\rho _{i}(e_{a},\xi _{i})\mathbb{B}(I_{k}e_{a},\xi _{k})+2\rho
_{k}(e_{a},\xi _{i})\mathbb{B}(I_{k}e_{a},\xi _{i})+T(\xi
_{i},e_{a},e_{b})(\nabla _{e_{b}}\mathbb{B})(I_{k}e_{a},\xi _{j}) \\
=4n(\nabla _{\xi _{i}}\mathbb{B})(\xi _{j},\xi _{k})-2(\nabla _{\xi
_{i}}L)(e_{a},e_{b})(L(I_{j}e_{a},I_{k}e_{b})+4n\zeta _{k}(\xi _{i},e_{a})%
\mathbb{B}(e_{a},\xi _{j}) \\
-2\rho _{i}(e_{a},\xi _{i})\mathbb{B}(I_{k}e_{a},\xi _{k})+2\rho
_{k}(e_{a},\xi _{i})\mathbb{B}(I_{k}e_{a},\xi _{i})+T(\xi
_{i},e_{a},e_{b})(\nabla _{e_{b}}\mathbb{B})(I_{k}e_{a},\xi _{j})
\end{multline}%
Substituting \eqref{int601} and \eqref{int602} in \eqref{int60} we
come to
\begin{multline}\label{e:c6.1}
4n\Big[(\nabla _{\xi _{i}}\mathbb{B})(\xi _{j},\xi _{k})-(\nabla _{\xi _{k}}%
\mathbb{B})(\xi _{i},\xi _{j})\Big]\\
= 2(\nabla _{e_{b}}L)(I_{k}e_{b},I_{j}e_{a})\mathbb{B}(e_{a},\xi
_{i}) + 2\Big[(\nabla _{e_{b}}\mathbb{B})(e_{a},\xi _{i})+(\nabla
_{\xi
_{i}}L)(e_{b},e_{a})\Big]L(I_{k}e_{b},I_{j}e_{a})\\-\mathbb{B}(e_{a},\xi
_{j}) \Big[4n\zeta _{k}(\xi _{i},e_{a})+(\nabla _{e_{b}}T)(\xi
_{i},I_{k}e_{b},e_{a})\Big]+ \Big[2\rho _{i}(e_{a},\xi
_{i})+\frac{1}{2n}d(tr\,L)(e_{a})\Big] \mathbb{B}(I_{k}e_{a},\xi
_{k})\\-2\rho _{k}(e_{a},\xi _{i})\mathbb{B}(I_{k}e_{a},\xi _{i})
+T(\xi_{i},I_{k}e_{a},e_{b})\Big[(\nabla
_{e_{b}}\mathbb{B})(e_{a},\xi _{j})- (\nabla
_{e_{a}}\mathbb{B})(e_{b},\xi _{j})\Big]
\\
+\frac{1}{2n}%
(tr\,L)(\nabla _{e_{a}}\mathbb{B})(I_{k}e_{a},\xi _{k})-4(tr\,L)\mathbb{B}%
(\xi _{j},\xi _{j})+4(tr\,L)\mathbb{B}(\xi _{i},\xi _{i})
\end{multline}%
 With the help
of \eqref{rverg}, the symmetry of $L$, and the divergence formulas
\eqref{e:div LI} and \eqref{e:div L} we find
\begin{multline}
4n\zeta _{k}(\xi _{i},e_{a})=(4n+1)\mathbb{B}(I_{k}e_{a},\xi
_{i})-\mathbb{B}(e_{a},\xi _{j})+\mathbb{B }(I_{i}e_{a},\xi
_{k})+\frac{1}{4n}d(trL)(I_{j}e_{a})\\
 =-\frac{1}{4}\left[ \nabla
L(e_{b},I_{i}e_{a},I_{k}e_{b})+\nabla
L(e_{b},e_{a},I_{j}e_{b})+\nabla
L(e_{b},I_{k}e_{a},I_{i}e_{b})+\nabla
L(e_{b},I_{j}e_{a},e_{b})\right]  .
\end{multline}
It follows from \eqref{t1} that
\begin{multline}\label{e:nablaTvhh}
(\nabla_{e_b}T)(\xi _{i},I_{k}e_{b},e_{a}) =\frac{1}{4}(\nabla_{e_b}L)(I_{j}{ e_{b}},{e_{a}})
-\frac{3}{4}(\nabla_{e_b}L)({I_{k}e_{b}},I_{i}{e_{a}})+\frac{1 }{4}(\nabla_{e_b}
L)({e_{b}},I_{j}e_{a})\\
+\frac{1}{4}(\nabla_{e_b}L)(I_{i}e_{b},I_{k}e_{a})
-\frac{1}{4n}d(trL)(I_{j}e_{a}).
\end{multline}
Adding the last two equations we see
\begin{equation*}%\label{comput3}
{4n\zeta }_{k}{(\xi }_{i}{,e}_{a}{)+(\nabla_{e_b}T)(\xi }_{i}{,I}_{k}{e}_{b}{,e}_{a})
=4nB(I_{k}{e}_{a}\xi_i)-4nB(I_{i}{e}_{a}{\xi}_{k}).
\end{equation*}
Using in \eqref{e:c6.1} the above identity, Lemma \ref{l:rho-ta},
\eqref{intexih11}, \eqref{inte1}, together with
$L(\mathbf{e}_{b},I_{s}e_{b})=0$, a long calculation gives
\begin{multline*}
4n(\nabla_{\xi_i}\mathbb B)(\xi _{k},\xi _{j})-4n(\nabla_{\xi_k}\mathbb B)(\xi
_{i},\xi _{j})=-4(trL)B(\xi _{j},\xi _{j})+2(trL)B(\xi _{i},\xi
_{i})+2(trL)B(\xi
_{k},\xi _{k})\\
8nB(I_{j}e_{a},\xi _{k})B(e_{a},\xi _{i})+4nB(e_{a},\xi
_{i})B(I_{k}e_{a},\xi _{j})+4nB(I_{i}e_{a},\xi _{k})B(e_{a},\xi
_{j}) .
\end{multline*}
Hence \eqref{intehxi41} is proven.

The other integrability condition in this case, \eqref{intehxi441},
can be obtained similarly using \eqref{intehxi312} and the Ricci
identities.
\end{proof}
The proof of Theorem~\ref{main2} is complete.

\section{A Ferrand-Obata type theorem}
%With the help
The group of conformal quaternionic contact automorphisms is a Lie
group, { which follows for example} from the equivalence of a qc
structure with a regular normal parabolic geometry. We thank A.
Cap for explaining this point to us.

{ A standard application of Theorem \ref{main1}, Theorem \ref{main2} { and Corollary \ref{main3}} (following \cite{W,W2}) gives a proof of
Ferrand-Obata type theorem concerning the { (4n+3) dimensional sphere} established in a more general situation for a parabolic structure
admitting { regular} Cartan connection in \cite{F}}.

Recall that the unit sphere $S$ in quaternionic space has a
natural  qc structure, namely, the standard 3-Sasakian
structure on the sphere}, cf. Section \ref{s:standard structures}.

\begin{thrm}\label{t:auto group}
Let $(M, \eta)$ be a compact quaternionic contact manifold and $G$ a
connected Lie group of conformal quaternionic contact automorphisms of
$M$. If $G$  is non-compact then $M$ is qc conformally equivalent to the
unit sphere $S$ in quaternionic space.
\end{thrm}

\begin{proof}
The argument { follows closely the proof of the CR case and} has two
steps. The first step is to show the local equivalence to the
sphere. This is done analogously to the CR case \cite{W} or
\cite{W2}. The key of the proof is the existence of an invariant
one form of contact structures when $M$ is not locally {
qc-conformally} flat. In our case this is achieved with the help
of the qc conformal curvature tensor $W^{qc}$, namely
$\eta^*=||W^{qc}||\eta$, cf. { Theorem \ref{main1} and
Corollary \ref{main3}}. Let $G_1$ be a one-parameter subgroup of
$G$ with infinitesimal generator $Q$. Suppose $M$ is not locally
flat, it is enough to show that $G_1$ is compact, which will be a
contradiction with the non-compactness of $G$ by \cite{MZ}.  Let
$U$ be a non-empty connected open set where $W^{qc}$ does not
vanish. Consider $\eta^*$ only on $U$, where it is an invariant
form under conformal qc transformations. Let
$\xi_1^*,\xi_2^*,\xi_3^*$ be the Reeb vector fields of the form
$\eta^*$. We have $\mathcal{L}_Q \eta_j^*=0$ and $\mathcal{L}_Q
\xi_j^*=0$. On $U$ we can decompose $Q$ uniquely as
$Q=Q_H+\eta^*_j(Q)\xi^*_j$. By \cite[Corollary 7.5]{IMV}, the
function $f=\sum_{j=1}^3[\eta^*_j(Q)]^2$ does not vanish
identically. For a sufficiently small $\epsilon$ define the set
$U_\epsilon=\{m\in U:\, f\geq \epsilon \}$. Notice that
$f(m)\rightarrow 0$ as $m\rightarrow\partial U$ and thus
$U_\epsilon$ is a closed subset of $M$ hence a compact.
Furthermore $U_\epsilon$ is invariant under the flow of $Q$ since
$Qf=\sum_{j=1}^3\eta^*_j(Q)[(\mathcal{L}_Q
\eta_j^*)(Q)+\eta^*_j([Q,Q])] = 0$, and also under the closure of
$G_1$. Let $P$ be the principle bundle over  $(U,\eta^*)$ with
fibre isomorphic to Sp(n)Sp(1) determined by the qc structure
$\eta^*$, and $P_\epsilon$ the part over $U_\epsilon$. From the
above considerations $P_\epsilon$ is invariant under the closure
$\bar G_1$. By \cite[Ch. I,Theorem 3.2]{Ko} $G_1$ embeds as a
closed submanifold of $P_\epsilon$. Hence $G_1$ is compact and the
proof of step 1 is complete.

The second step is to show the global equivalence, which is done
in \cite[Proposition D]{K}, see also \cite[Problem E (a)]{K} and
for the CR case \cite{L}. The analysis {there} involves the
dynamics of one-parameter groups of qc conformal automorphisms.
\end{proof}

\begin{rmrk}
a) The dynamics of one-parameter groups of conformal automorphisms
has been studied in the more general setting of boundaries of rank
one symmetric spaces in \cite{F}, which results in the proof of a
general Ferrand-Obata type theorem { established in \cite{F}}.

b) We note that the  properties of the curvature of the Biquard
connection investigated in the paper allow to apply the analysis of
\cite{S} to the case of  a qc structure { (Proposition~\ref{hflat}
supports \cite[Lemma~4.1]{S}, the quaternionic contact parabolic
normal coordinates of \cite{Ku} supply the pseudohermitian normal
coordinates of \cite{JL} in the quaternionic contact  case) and,
thus,  a proof of Theorem~\ref{t:auto group} could be obtained
following the approach of \cite{S}}.

\end{rmrk}

\end{document}